\pdfoutput=1
\RequirePackage{ifpdf}
\ifpdf 
\documentclass[pdftex]{sigma}
\else
\documentclass{sigma}
\fi

\usepackage{enumerate}

\numberwithin{equation}{section}

\newtheorem{Theorem}{Theorem}[section]
\newtheorem{Lemma}[Theorem]{Lemma}
\newtheorem{Proposition}[Theorem]{Proposition}
 { \theoremstyle{definition}
\newtheorem{Definition}[Theorem]{Definition}
\newtheorem{Example}[Theorem]{Example}
\newtheorem{Remark}[Theorem]{Remark} }

\DeclareMathOperator{\Ima}{Im}

\begin{document}

\newcommand{\arXivNumber}{1409.8177}

\allowdisplaybreaks

\renewcommand{\thefootnote}{$\star$}

\renewcommand{\PaperNumber}{049}

\FirstPageHeading

\ShortArticleName{A Combinatorial Formula for Certain Elements of Upper Cluster Algebras}

\ArticleName{A Combinatorial Formula for Certain Elements \\ of Upper Cluster Algebras\footnote{This paper is a~contribution to the Special Issue on New Directions in Lie Theory.
The full collection is available at \href{http://www.emis.de/journals/SIGMA/LieTheory2014.html}{http://www.emis.de/journals/SIGMA/LieTheory2014.html}}}

\Author{Kyungyong LEE~$^{\dag\ddag}$, Li LI~$^\S$ and Matthew R. MILLS~$^{\dag}$}

\AuthorNameForHeading{K.~Lee, L.~Li and M.R.~Mills}

\Address{$^\dag$~Department of Mathematics, Wayne State University, Detroit, MI 48202, USA}
\EmailD{\href{mailto:klee@math.wayne.edu}{klee@math.wayne.edu}, \href{mailto:matthew.mills2@wayne.edu}{matthew.mills2@wayne.edu}}

\Address{$^\ddag$~Korea Institute for Advanced Study, Seoul, Republic of Korea 130-722}
\EmailD{\href{mailto:klee1@kias.re.kr}{klee1@kias.re.kr}}

\Address{$^\S$~Department of Mathematics and Statistics,
Oakland University,
Rochester, MI 48309, USA}
\EmailD{\href{mailto:li2345@oakland.edu}{li2345@oakland.edu}}

\ArticleDates{Received September 30, 2014, in f\/inal form June 22, 2015; Published online June 26, 2015}

\Abstract{We develop an elementary formula for certain non-trivial elements  of upper cluster algebras. These elements have positive coef\/f\/icients. We show that when the cluster algebra is acyclic these elements form a basis. Using this formula, we show that each non-acyclic skew-symmetric cluster algebra of rank 3 is properly contained in its upper cluster algebra. }


\Classification{13F60}

\renewcommand{\thefootnote}{\arabic{footnote}}
\setcounter{footnote}{0}

\section{Introduction}
Cluster algebras were introduced by Fomin and Zelevinsky in \cite{FZ1}. A cluster algebra $\mathcal{A}$ is a~subalgebra of a rational function f\/ield with a distinguished set of generators, called \emph{cluster variables}, that are generated by an iterative procedure called \emph{mutation}. By construction cluster variables are rational functions, but it is shown in \emph{loc.~cit.}~that  they are Laurent
polynomials with integer coef\/f\/icients. Moreover, these coef\/f\/icients are known to be non-negative \cite{GHKK, LS}.

\looseness=-1
Each cluster algebra $\mathcal{A}$ also determines an \emph{upper cluster algebra} $\mathcal{U}$, where $\mathcal{A}\subseteq \mathcal{U}$ \cite{BFZ}. It is believed, especially in the context of algebraic geometry, that $\mathcal{U}$ is better behaved than  $\mathcal{A}$ (for instance, see \cite{BMRS, GHK, GHKK}). Matherne and Muller \cite{MM} gave a general algorithm to compute generators of $\mathcal{U}$. Plamondon \cite{P1,P2} obtained a (not-necessarily positive) formula for certain elements of skew-symmetric upper cluster algebras using quiver representations.    However a directly computable and manifestly positive formula for (non-trivial) elements in $\mathcal{U}$ is not available yet.

In this paper we develop an elementary formula for a family of elements $\{\tilde{x}[{\bf a}]\}_{{\bf a}\in\mathbb{Z}^n}$ of the upper cluster algebra for any f\/ixed initial seed $\Sigma$. We write $\tilde{x}_\Sigma[{\bf a}]$ for $\tilde{x}[{\bf a}]$ when we need to emphasize the dependence on the initial seed $\Sigma$.
This family of elements are constructed in Def\/inition~\ref{definition:GCC2} in terms of sequences of sequences, and an equivalent def\/inition is given in terms of Dyck paths and globally compatible collections in Def\/inition~\ref{definition:GCC}.
One of our main theorems is the following. For the def\/inition of {\it geometric type}, see Section~\ref{section2}.

\begin{Theorem}\label{thm:x[a]}
Let $\mathcal{U}$ be the upper cluster algebra of a $($not necessarily acyclic$)$ cluster algebra~$\mathcal{A}$ of geometric type, and~$\Sigma$ be any seed. Then  $\tilde{x}_\Sigma[\mathbf{a}]\in\mathcal{U}$ for all~$\mathbf{a}\in \mathbb{Z}^n$.
\end{Theorem}

These elements have some nice properties. They have positive coef\/f\/icients by def\/inition; they are multiplicative in the sense that we can factorize an element  $\tilde{x}[{\bf a}]$ (${\bf a}\in\mathbb{Z}_{\ge0}^n$) into ``elementary pieces''  $\tilde{x}[{\bf a}']$ where all entries of ${\bf a}'$ are $0$ or $1$;
for an equioriented quiver of type $A$, these elements form a canonical basis~\cite{B9}. Moreover, we shall prove in Section~\ref{section6} the following result. (For the terminology and notation therein, see Section~\ref{section2}.)

\begin{Theorem}\label{thm:basis}
Let $\mathcal{A}$ be an acyclic cluster algebra of geometric type, and $\Sigma$ be an acyclic seed. Then $\{\tilde{x}_\Sigma[\mathbf{a}]\}_{\mathbf{a}\in\mathbb{Z}^n}$  form a $\mathbb{ZP}$-basis of $\mathcal{A}$.
\end{Theorem}

For a non-acyclic seed $\Sigma$, the family $\{\tilde{x}_\Sigma[\mathbf{a}]\}$ may neither span $\mathcal{U}$ nor be linearly independent. For a linearly dependent example, see Example~\ref{example}(b). Nevertheless, for certain non-acyclic cluster algebras and for some choice of ${\bf a} \in \mathbb{Z}^n$, the element $\tilde{x}_\Sigma[{\bf a}]$ can be used to construct elements in $\mathcal{U}\setminus \mathcal{A}$. One of our main results in this direction is the following:

\begin{Theorem}\label{main-cor}
A non-acyclic rank three skew-symmetric cluster algebra $\mathcal{A}$ of geometric type is not equal to its upper cluster algebra $\mathcal{U}$.
\end{Theorem}

This theorem  is inspired by the following results: Berenstein, Fomin and Zelevinsky  \cite[Proposition 1.26]{BFZ} showed that $\mathcal{A}\neq \mathcal{U}$ for the Markov skew-symmetric matrix $M(2)$, where
\begin{gather*}
M(a)=\left(\begin{matrix}0 & a & -a\\ -a & 0 & a\\ a & -a & 0\end{matrix}\right).
\end{gather*}
Speyer \cite{S} found an inf\/initely generated upper cluster algebra, which is the one associated to the skew-symmetric matrix $M(3)$ with generic coef\/f\/icients. On the contrary, \cite[Proposition~6.2.2]{MM}  showed that the upper cluster algebra associated to $M(a)$ for $a\geq 2$ but with trivial coef\/f\/icients is f\/initely generated, which implies $\mathcal{A}\neq \mathcal{U}$ because $\mathcal{A}$ is known to be inf\/initely generated \cite[Theorem~1.24]{BFZ}.

The paper is organized as follows. In Section~\ref{section2} we review def\/initions of cluster algebras and upper cluster algebras. Section~\ref{section3} is devoted to the construction of $\tilde{x}[{\bf a}]$ and the proof of Theorem~\ref{thm:x[a]}, and Section~\ref{section4} to the proof of
Theorem \ref{main-cor}, that $\mathcal{A}\neq \mathcal{U}$ for non-acyclic rank~3 skew-symmetric cluster algebras. Section~\ref{section5} introduces the Dyck path formula and its relation with the construction in Section~\ref{section3}. Finally, in Section~\ref{section6}, we present two proofs of Theorem~\ref{thm:basis}.

\section{Cluster algebras and upper cluster algebras}\label{section2}

Let $m$, $n$ be positive integers such that $m\geq n$. Denote $\mathcal{F}=\mathbb{Q}(x_1,\dots,x_m)$.
A \emph{seed}  $\Sigma=(\tilde{\bf x},\tilde{B})$ is a pair where $\tilde{\bf x}=\{x_1,\ldots,x_m\}$ is an $m$-tuple of elements of $\mathcal{F}$ that form a free generating set, $\tilde{B}$ is an $m\times n$ integer matrix such that the submatrix $B$  (called the \emph{principal part}) formed by the top $n$ rows is sign-skew-symmetric (that is, either $b_{ij}=b_{ji}=0$, or else~$b_{ji}$ and~$b_{ij} $ are of opposite sign; in particular, $b_{ii}=0$ for all~$i$).
The integer $n$ is called the \emph{rank} of the seed.

For any integer $a$, let $[a]_+:= \max(0,a)$.
 Given a seed $(\tilde{\bf x},\tilde{B})$ and a specif\/ied index $1\leq k \leq n$, we def\/ine \emph{mutation} of $(\tilde{\bf x},\tilde{B})$ at $k$, denoted $\mu_{k}(\tilde{\bf x},\tilde{B})$, to be a new seed $(\tilde{\bf x}',\tilde{B}')$, where
\begin{gather*}
 x_i' = \begin{cases}
\displaystyle  x'_k=x_k^{-1}\left(\prod\limits_{i=1}^{m}x_i^{[b_{ik}]_+}+\prod\limits_{i=1}^{m}x_i^{[-b_{ik}]_+}\right),
  & \mbox{if} \ \  i=k, \\
x_i, & \mbox{otherwise}, \end{cases}\\
b_{ij}' = \begin{cases} -b_{ij}, & \mbox{if} \ \ i=k \ \ \text{or} \ \ j=k,  \\
b_{ij} + \dfrac{|b_{ik}|b_{kj} + b_{ik}|b_{kj}|}{2}, & \mbox{otherwise.}
\end{cases}
\end{gather*}
If the principal part of $\tilde{B}'$ is also sign-skew-symmetric, we say that the mutation is well-def\/ined.
Note that a well-def\/ined mutation is an involution, that is mutating $(\tilde{\bf x}', \tilde{B}')$ at~$k$ will return to our original seed $(\tilde{\bf x}, \tilde{B})$.

Two seeds $\Sigma_1$ and $\Sigma_2$ are said to be \emph{mutation-equivalent} or in the same \emph{mutation class} if~$\Sigma_2$ can be obtained by a sequence of well-def\/ined mutations from~$\Sigma_1$. This is obviously an equivalence relation. A seed~$\Sigma$ is said to be \emph{totally mutable}  if every sequence of mutations from~$\Sigma$ consists of well-def\/ined ones. It is shown in \cite[Proposition~4.5]{FZ1} that a seed is totally mutable if $B$ is skew-symmetrizable, that is, if there exists a diagonal matrix $D$ with positive diagonal entries such that $DB$ is skew-symmetric.

To emphasize the dif\/ferent roles played by $x_i$ $(i\le n)$ and $x_i$ $(i>n)$, we also use $({\bf x},{\bf y},B)$ to denote the seed $(\tilde{\bf x}, \tilde{B})$, where ${\bf x}=\{x_1,\dots,x_n\}$, ${\bf y}=\{y_1,\dots,y_n\}$ where $y_j=\prod\limits_{i=n+1}^mx_i^{b_{ij}}$. We call~${\bf x}$ a~\emph{cluster}, ${\bf y}$ a~\emph{coefficient tuple}, $B$ the \emph{exchange matrix}, and the elements of a cluster \emph{cluster variables}. We denote
\begin{gather*}
\mathbb{ZP}=\mathbb{Z}\big[x_{n+1}^{\pm1},\dots,x_m^{\pm1}\big].
\end{gather*}

In the paper, we shall only study cluster algebras of geometric type, def\/ined as follows.

\begin{Definition}
Given a totally mutable seed $({\bf x}, {\bf y},B)$, the \emph{cluster algebra} $\mathcal{A}({\bf x},{\bf y},B)$ \emph{of geometric type} is the subring of $\mathcal{F}$ generated over $\mathbb{ZP}$ by
\begin{gather*}
\bigcup_{({\bf x}',{\bf y}',B')} {\bf x}',
 \end{gather*}
 where the union runs over all seeds $({\bf x}', {\bf y}',B')$ that are mutation-equivalent to $({\bf x},{\bf y},B)$. The seed $({\bf x},{\bf y},B)$ is called the \emph{initial seed} of $\mathcal{A}({\bf x},{\bf y},B)$. (Since $({\bf x},{\bf y},B)=(\tilde{\bf x},\tilde{B})$ in our notation, $\mathcal{A}({\bf x},{\bf y},B)$ is also denoted $\mathcal{A}(\tilde{\bf x},\tilde{B})$.)
\end{Definition}

It follows from the def\/inition that any seed in the same mutation class will generate the same cluster algebra up to isomorphism.

 For any $n \times n$ sign-skew-symmetric matrix $B$, we associate a (simple) directed graph $Q_B$ with vertices $1,\ldots,n$, such that for each pair $(i,j)$ with $b_{ij}>0$, there is  \emph{exactly one} arrow from vertex~$i$ to vertex~$j$.
(Note that even if $B$ is skew-symmetric, $Q_B$ is not the usual quiver associated to~$B$ which can have multiple edges.)

 We call $B$ (as well as the digraph $Q_B$ and the seed $\Sigma=({\bf x},{\bf y},B)$) \emph{acyclic} if there are no oriented cycles in $Q_B$.
We say that the cluster algebra $\mathcal{A}({\bf x},{\bf y},B)$ is \emph{acyclic} if there exists an acyclic seed; otherwise we say that the cluster algebra is \emph{non-acyclic}.

\begin{Definition}
Given a cluster algebra $\mathcal{A}$, the \emph{upper cluster algebra} $\mathcal{U}$ is def\/ined as
\begin{gather*}
 \mathcal{U} =
 \bigcap_{{\bf x}=\{x_1,\ldots,x_n\}}
  \mathbb{ZP}\big[x_1^{\pm 1},\ldots,x_n^{\pm 1}\big],
\end{gather*}
where ${\bf x}$ runs over all clusters of~$\mathcal{A}$.
\end{Definition}

Now we can give the following def\/inition of coprime when the cluster algebra is of geometric type (given in \cite[Lemma~3.1]{BFZ}).

\begin{Definition}
A seed  $(\tilde{\bf x},\tilde{B})$ is \emph{coprime} if no two columns of $\tilde{B}$ are proportional to each other with the proportionality coef\/f\/icient being a ratio of two odd integers.

A cluster algebra is \emph{totally coprime} if every seed is coprime.
\end{Definition}

In certain cases it is suf\/f\/icient to consider only the clusters of the initial seed and the seeds that are a single mutation away from it, rather than all the seeds in the entire mutation class.
For a cluster ${\bf x}$,  let $\mathcal{U}_{{\bf x}}$ be the intersection in $\mathbb{ZP}(x_1,\ldots,x_n)$ of the $n+1$ Laurent rings corresponding to ${\bf x}$ and its one-step mutations:
\begin{gather*}
\mathcal{U}_{\bf x} := \mathbb{ZP}\big[x_1^{\pm 1},\ldots, x_n^{\pm 1}\big] \cap \left(\bigcap_i \mathbb{ZP}\big[x_1^{\pm 1},\ldots,x_i'^{\pm 1},\ldots, x_n^{\pm 1}\big]\right).
\end{gather*}

\begin{Theorem}[\cite{BFZ,M}]\label{UB}
We have
 $\mathcal{A}\subseteq\mathcal{U}\subseteq\mathcal{U}_{\bf x}$. Moreover,
\begin{enumerate}\itemsep=0pt
\item[$(i)$] If $\mathcal{A}$ is acyclic, then $\mathcal{A}=\mathcal{U}$.

\item[$(ii)$] If $\mathcal{A}$ is totally coprime, then $\mathcal{U}=\mathcal{U}_{{\bf x}}$ for any seed $({\bf x},{\bf y},B)$. In particular, this holds when the matrix~$\tilde{B}$ has full rank.
\end{enumerate}
\end{Theorem}

\section{Construction of some elements in the upper cluster algebra}\label{section3}
Fix an initial seed $\Sigma$ (thus $\tilde{B}$ is f\/ixed). In this section, we construct Laurent polyno\-mials~$\tilde{x}[{\bf a}]$ $(=\tilde{x}_\Sigma[{\bf a}])$ and show that they are in the upper cluster algebra.

We def\/ine $b_{ij}=-b_{ji}$ for $1\le i\le n$, $n+1\le j\le m$, and def\/ine $b_{ij}=0$ if $i,j>n$.
Def\/ine
\begin{gather*}
Q_{\tilde B}=\{(i,j)\, |\, 1\le i,j\le m,\;  b_{ij}>0\}.
\end{gather*}
By abuse of notation we also use~$Q_{\tilde B}$ to denote the digraph with vertex set~$\{1,\dots,m\}$ and edge set~$Q_{\tilde B}$. Then~$Q_B$ is a full sub-digraph of~$Q_{\tilde B}$ that consists of the f\/irst $n$ vertices.

We def\/ine the following operations on the set of f\/inite $\{0,1\}$-sequences. Let $t=(t_1,\dots,t_a)$, $t'=(t'_1,\dots,t'_b)$. Def\/ine
\begin{gather}\label{eq:t sum and dot}
\bar{t}=(\bar{t}_1,\dots,\bar{t}_a)=(1-t_1,\dots,1-t_a),\qquad |t|=\sum_{r=1}^a t_r,\qquad t\cdot t'=\sum_{r=1}^{\min(a,b)} t_rt'_r.
\end{gather}
and for $t=()\in \{0,1\}^0$, def\/ine $\bar{t}=()$.

\begin{Definition}\label{definition:GCC2}
Let $\mathbf{a}=(a_i)\in \mathbb{Z}^n$.
\begin{enumerate}\itemsep=0pt
\item[(i)] Let $\mathbf{s}=(\mathbf{s}_1,\dots,\mathbf{s}_n)$ where
$\mathbf{s}_i=(s_{i,1},s_{i,2},\dots,s_{i,[a_i]_+})\in\{0,1\}^{[a_i]_+}$ for $i=1,\dots,n$. Let $S_{\rm all}=S_{\rm all}({\bf a})$ be the set of all such $\mathbf{s}$. Let $S_{\rm gcc}=S_{\rm gcc}({\bf a})=\{\mathbf{s}\in S_{\rm all}\, |\, \mathbf{s}_i\cdot\bar{\mathbf{s}}_j=0$ for every $(i,j)\in Q_B\}$.

By convention, we assume that $a_i=0$ and $\mathbf{s}_i=()$ for $i>n$.

\item[(ii)] Def\/ine $\tilde{x}[\mathbf{a}]$ $(=\tilde{x}_\Sigma[{\bf a}])$ to be the Laurent polynomial
\begin{gather*}
\tilde{x}[\mathbf{a}]:=\left(\prod_{l=1}^n  x_l^{-a_l}\right)
\sum_{\mathbf{s}\in S_{\rm gcc}}\left(\prod_{(i,j)\in Q_{\tilde B}}  x_{i}^{b_{ij}| \bar{\mathbf{s}}_j |} x_{j}^{-b_{ji}|\mathbf{s}_i|}\right).
\end{gather*}
(Note that the exponent $-b_{ji}|\mathbf{s}_i|$ is nonnegative since $b_{ji}<0$ for $(i,j)\in Q_{\tilde B}$.)

\item[(iii)] Def\/ine $z[\mathbf{a}]$ $(=z_\Sigma[{\bf a}])$ to be the Laurent polynomial
\begin{gather*}
z[\mathbf{a}]:=\left(\prod_{l=1}^n  x_l^{-a_l}\right)
\sum_{\mathbf{s}\in S_{\rm all}}\left(\prod_{(i,j)\in Q_{\tilde B}}  x_{i}^{b_{ij}|\bar{\mathbf{s}}_j|} x_{j}^{-b_{ji}|\mathbf{s}_i|}\right).
\end{gather*}
\end{enumerate}
\end{Definition}

\begin{Example}\label{example}
(a) Use Def\/inition \ref{definition:GCC2} to compute $\tilde{x}[\mathbf{a}]$ for $n=2$, $m=3$, $\mathbf{a}=(1,1)$, and
\begin{gather*}
\tilde{B}=\begin{bmatrix}
0&a\\
-a'&0\\
c&-b
\end{bmatrix},\qquad a,a',b,c>0.
\end{gather*}
Then we have $Q_{\tilde B}=\{(1,2),(2,3),(3,1)\}$, $S_{\rm gcc}=\{((0),(0)), ((0),(1)),  ((1),(1))\}$. Note that $\mathbf{s}=((1),(0))$ is not in $S_{\rm gcc}$ because $(1,2)\in Q_{\tilde B}$ but $\mathbf{s}_1\cdot\bar{\mathbf{s}}_2=1(1-0)=1\neq 0$. Thus
\begin{gather*}
\tilde{x}[\mathbf{a}]
 =x_1^{-1}x_2^{-1}
\big( \big(x_1^ax_2^0\big)\big(x_2^0x_3^0\big)\big(x_3^cx_1^0\big)+
 \big(x_1^0x_2^0\big)\big(x_2^0x_3^b\big)\big(x_3^cx_1^0\big)+
  \big(x_1^0x_2^{a'}\big)\big(x_2^0x_3^b\big)\big(x_3^0x_1^0\big)\big)\\
\hphantom{\tilde{x}[\mathbf{a}]}{}
=x_1^{-1}x_2^{-1}
\big(x_1^ax_3^c+x_3^bx_3^c+x_2^{a'} x_3^b\big).
\end{gather*}

(b) For a non-acyclic seed, $\tilde{x}[\mathbf{a}]$ is less interesting for certain choices of ${\bf a}$: take $n=m=3$, $\mathbf{a}=(1,1,1)$ and
\begin{gather*}
\tilde{B}=\begin{bmatrix}
0&a&-c'\\
-a'&0&b\\
c&-b'&0
\end{bmatrix},\qquad a,a',b,b',c,c'>0.
\end{gather*}
Then $Q_{\tilde B}=\{(1,2),(2,3),(3,1)\}$, $S_{\rm gcc}=\{((0),(0),(0)), ((1),(1),(1))\}$. Thus
\begin{gather*}
\tilde{x}[\mathbf{a}]
=\frac{x_1^ax_2^bx_3^c+x_2^{a'}x_3^{b'}x_1^{c'}}{x_1x_2x_3},
\end{gather*}
which can be reduced to $x_1^{a-1}x_2^{b-1}x_3^{c-1}+x_2^{a'-1}x_3^{b'-1}x_1^{c'-1}$, that is
\begin{gather*}
\tilde{x}[1,1,1]=\tilde{x}[(1-a,1-b,1-c)]+\tilde{x}[(1-a',1-b',1-c')].
\end{gather*}
Thus $\{\tilde{x}[{\bf a}]\}$ is not $\mathbb{ZP}$-linear independent.
\end{Example}

\begin{Lemma}[multiplicative property of $\tilde{x}{[{\bf a}]}$ and $z{[{\bf a}]}$]\label{multiplicative} Fix a seed $\Sigma$.
\begin{enumerate}\itemsep=0pt
\item[$(i)$] For $k,a\in\mathbb{Z}$, ${\bf a}\in\mathbb{Z}^n$, define
\begin{gather*}
f_k(a)
:=\begin{cases}
1,&\textrm{if} \ \ k\le a,\\
0, &\textrm{otherwise},
\end{cases}
\end{gather*}
$f_k(\mathbf{a}):=(f_k(a_1),\dots,f_k(a_n))\in\{0,1\}^n$,
$\mathbf{a}_+:=([a_1]_+, [a_2]_+,\dots, [a_n]_+)$.
Then
\begin{gather*}
\tilde{x}[\mathbf{a}]
=\left(\prod_{i=1}^n  x_i^{[-a_i]_+}\right)\tilde{x}[\mathbf{a}_+]
=\left(\prod_{i=1}^n  x_i^{[-a_i]_+}\right)\prod_{k\ge1}
\tilde{x}[f_k(\mathbf{a}_+)].
\end{gather*}

\item[$(ii)$] Assume that the underlying undirected graph of $Q_B$ has $c>1$ components, inducing a~partition of the vertex set $\{1,\dots,n\}=I_1\cup\cdots\cup I_c$. Define
${\bf a}^{(j)}=\big(a^{(j)}_1,\dots,a^{(j)}_c\big)$ by $a^{(j)}_i=a_i$ if $i\in I_j$, otherwise $a^{(j)}_i=0$. Then
\begin{gather*}
\tilde{x}[{\bf a}]=\prod_{j=1}^c\tilde{x}\big[{\bf a}^{(j)}\big].
\end{gather*}
$($Note that each factor $\tilde{x}[{\bf a}^{(j)}]$ can be regarded as an element in a cluster algebra of rank $|I_j|<n$, with the same $\mathbb{ZP}.)$

\item[$(iii)$]
We have $z[{\bf a}]=\prod\limits_{i=1}^nx_i^{\langle-a_i\rangle}$ where we use the notation
\begin{gather*}
x_i^{\langle r\rangle}=\begin{cases}
x_i^r,  &\textrm{if} \ \ r\ge 0,\\
(x_i')^{-r}, &\textrm{if} \ \ r< 0.
\end{cases}
\end{gather*}
$($Recall that $x_i'$ is obtained by mutating the initial seed at $i)$.
As a consequence, $(i)$, $(ii)$ still hold if we replace~$\tilde{x}[-]$ by~$z[-]$. Moreover,
if the seed $\Sigma$ is acyclic, then $\{z[\mathbf{a}]\}_{\mathbf{a}\in\mathbb{Z}^n}$ is the standard monomial basis, i.e., the set of monomials in $x_1,\dots,x_n$, $x_1',\dots,x_n'$  which contain no product of the form~$x_jx_j'$.
\end{enumerate}
 \end{Lemma}

\begin{proof} (i)~The f\/irst equality is obvious. For the second equality, assuming ${\bf a}={\bf a}_+$.
For a~sequence $t=(t_1,\dots,t_a)\in\{0,1\}^a$, we can regard it as an inf\/inite sequence $(t_1,\dots,t_a,0,0,\dots)$. Then $\bar{t}=(f_r(a)-t_r)_{r=1}^\infty$. The sum and dot product in~\eqref{eq:t sum and dot} extend naturally.

Then
$\tilde{x}[\mathbf{a}]$ is the coef\/f\/icient of $z^0$ in the following polynomial in $\mathbb{Z}[x_1^{\pm 1},\ldots,x_m^{\pm 1}][z]$ (note that all the products $\prod\limits_{k\ge 1}$ appearing below are f\/inite products since the factors are $1$ if $k>\max([a_1]_+,\dots,[a_n]_+)$)
\begin{gather*}
\left(\prod_{i=1}^n  x_i^{-a_i}\right)
\sum_{\mathbf{s}\in S_{\rm all}}\left(\prod_{(i,j)\in Q_{\tilde B}}  x_{i}^{b_{ij}|\bar{\mathbf{s}}_j|} x_{j}^{-b_{ji}|\mathbf{s}_i|}z^{\mathbf{s}_i\cdot \bar{\mathbf{s}}_j}\right)\\
\qquad {} =\left(\prod_{i=1}^n  x_i^{-a_i}\right)
\sum_{\mathbf{s}\in S_{\rm all}}\prod_{k\ge1}\left(\prod_{(i,j)\in Q_{\tilde B}}  x_{i}^{b_{ij}\bar{s}_{j,k}} x_{j}^{-b_{ji}s_{i,k}}z^{s_{i,k}\bar{s}_{j,k}}\right)\\
\qquad {} =\left(\prod_{i=1}^n  x_i^{-a_i}\right)
\sum_{\mathbf{s}\in S_{\rm all}}\prod_{k\ge1}\left(\prod_{(i,j)\in Q_{\tilde B}}  x_{i}^{b_{ij}(f_k(a_j)-s_{j,k})} x_{j}^{-b_{ji}s_{i,k}}z^{s_{i,k}(f_k(a_j)-s_{j,k})}\right)\\
\qquad {} =\left(\prod_{i=1}^n  x_i^{-a_i}\right)\prod_{k\ge1}
\left(\sum_{\mathbf{s}^k\in S_{\rm all}^k}\prod_{(i,j)\in Q_{\tilde B}}  x_{i}^{b_{ij}(f_k(a_j)-s_{j,k})} x_{j}^{-b_{ji}s_{i,k}}z^{s_{i,k}(f_k(a_j)-s_{j,k})}\right),
\end{gather*}
where $S^k_{\rm all}$ is the set of all possible $\mathbf{s}^k=(s_{1,k},\dots,s_{n,k})$ with $(\mathbf{s}_1,\dots,\mathbf{s}_m)$ running through $S_{\rm all}$ (recall that~$s_{i,k}$ is the $k$-th number of ${\bf s}_i$, and by convention $s_{i,k}=0$ if $k>[a_i]_+$). Equivalently,
\begin{gather*}
S^k_{\rm all}=\big\{\mathbf{s}^k=(s_{1,k},\dots,s_{n,k})\in\{0,1\}^n\,|\,0\le s_{i,k}\le f_k(a_i) \; \textrm{for} \; i=1,\dots,n\big\}.
\end{gather*}
Meanwhile, denote $f_k(\mathbf{a})=(f_k(a_1),\dots,f_k(a_n))\in\{0,1\}^n$. Then $\tilde{x}[f_k(\mathbf{a})]$ is the coef\/f\/icient of~$z^0$ of
\begin{gather*}
\left(\prod_{i=1}^n  x_i^{-f_k(a_i)}\right)
\sum_{\mathbf{s}^k\in S_{\rm all}^k}\left(\prod_{(i,j)\in Q_{\tilde B}}  x_{i}^{b_{ij}(f_k(a_j)-s_{j,k})} x_{j}^{-b_{ji}s_{i,k}}z^{s_{i,k}(f_k(a_j)-s_{j,k})}\right).
\end{gather*}
So we conclude that
\begin{gather*}
\tilde{x}[\mathbf{a}]
 =\left(\prod_{i=1}^n x_i^{-a_i}\right)\prod_{k\ge1}\left(\tilde{x}[f_k(\mathbf{a})]\prod_{i=1}^n x_i^{f_k(a_i)}\right)\\
\hphantom{\tilde{x}[\mathbf{a}]}{} =\left(\prod_{i=1}^n x_i^{-a_i+\sum\limits_{k\ge1}f_k(a_i)}\right)\prod_{k\ge1}
\tilde{x}[f_k(\mathbf{a})]
=\prod_{k\ge1}
\tilde{x}[f_k(\mathbf{a})].
\end{gather*}

(ii) There is a bijection
\begin{gather*}
\prod_{j=1}^c S_{\rm gcc}\big({\bf a}^{(j)}\big)\to S_{\rm gcc}({\bf a}), \qquad
\big({\bf s}^{(1)},\dots, {\bf s}^{(c)}\big)\mapsto {\bf s}=({\bf s}_1,\dots,{\bf s}_n),
\end{gather*}
where ${\bf s}_i={\bf s}^{(j)}_i$ if $i\in I_j$.
This bijection induces the expected equality.

(iii) Rewrite
\begin{gather*}
z[{\bf a}] =\left(\prod_{i=1}^n  x_i^{-a_i}\right)
\sum_{{\bf s}_1,\dots,{\bf s}_n}\left(\prod_{i,j}  x_{i}^{|\bar{\mathbf{s}}_j|[b_{ij}]_+} x_{j}^{|\mathbf{s}_i|[-b_{ji}]_+}\right)\\
\hphantom{z[{\bf a}]}{} =
\sum_{{\bf s}_1,\dots,{\bf s}_n}\prod_{k=1}^n
\left(x_k^{-a_k}\prod_{i}  x_{i}^{|\bar{\mathbf{s}}_k|[b_{ik}]_+}\prod_j x_{j}^{|\mathbf{s}_k|[-b_{jk}]_+}\right)\\
\hphantom{z[{\bf a}]}{}
=
\prod_{k=1}^n
x_k^{-a_k}\sum_{{\bf s}_k}
\left(\prod_{i}  x_{i}^{|\bar{\mathbf{s}}_k|[b_{ik}]_+}\prod_j x_{j}^{|\mathbf{s}_k|[-b_{jk}]_+}\right)=\prod_{k=1}^n z[a_ke_k].\\
\end{gather*}
If $a_k\le 0$, then ${\bf s}_k=()$, therefore $z[a_ke_k]=x_k^{-a_k}=x_k^{\langle -a_k\rangle}$. If $a_k> 0$, then
\begin{gather*}
z[a_ke_k]
 =
x_k^{-a_k}\sum_{s_{k,1},\dots,s_{k,a_k}}
\left(\prod_{i}  x_{i}^{\sum\limits_{r=1}^{a_k}(1-s_{k,r})[b_{ik}]_+}\prod_j x_{j}^{\sum\limits_{r=1}^{a_k}(s_{k,r})[-b_{jk}]_+}\right)\\
\hphantom{z[a_ke_k]}{}
 =
x_k^{-a_k}
\prod_{r=1}^{a_k}
\sum_{s_{k,r}\in\{0,1\}}
\left(\prod_{i}  x_{i}^{(1-s_{k,r})[b_{ik}]_+}\prod_j x_{j}^{(s_{k,r})[-b_{jk}]_+}\right)\\
\hphantom{z[a_ke_k]}{}  =
x_k^{-a_k}
\prod_{r=1}^{a_k}
\left(\prod_{i}  x_{i}^{[b_{ik}]_+}+\prod_j x_{j}^{[-b_{jk}]_+}\right)
=(x'_k)^{a_k}=x_k^{\langle -a_k\rangle}.
\end{gather*}
This proves $z[{\bf a}]=\prod\limits_{i=1}^nx_i^{\langle-a_i\rangle}$. The analogue of~(i),~(ii) immediately follows. The fact that $\{z[\mathbf{a}]\}_{\mathbf{a}\in\mathbb{Z}^n}$ forms a basis is proved in  \cite[Theorem~1.16]{BFZ}.
\end{proof}

\begin{Remark}
For readers who are familiar with \cite[Lemma~4.2]{B9}, they may notice that the decomposition therein is similar to Lemma~\ref{multiplicative}(i) above. Indeed, \cite[Lemma~4.2]{B9} gives a~f\/iner decomposition. For example, for the coef\/f\/icient-free cluster algebra of the quiver $1\to 2\to 3$, $\tilde{x}[(2,1,3)]$ will decompose as $\tilde{x}[(1,1,1)]\tilde{x}[(1,0,1)]\tilde{x}[(0,0,1)]$ in  Lemma~\ref{multiplicative}(i),  but decompose as $\tilde{x}[(1,1,1)]\tilde{x}[(1,0,0)]\tilde{x}[(0,0,1)]^2$ in \cite[Lemma~4.2]{B9}.
\end{Remark}

The following lemma focuses on the case where $\mathbf{a}=(a_i)\in \{0,1\}^n$ as opposed to that in~$\mathbb{Z}^n$. The condition of~$S_{\rm gcc}$ takes a much simpler form: we can treat sequences~$(0)$ and~$(1)$ as numbers~$0$ and~$1$ respectively, and the condition $\mathbf{s}_i\cdot\bar{\mathbf{s}}_j=0$ can be written as $(s_i,a_j-s_j)\neq(1,1)$. We shall use $\mathbf{S}$ to denote this simpler form of~$S_{\rm gcc}$.

\begin{Lemma}\label{lemma:01}
For $\mathbf{a}=(a_i)\in \{0,1\}^n$, denote by
$\mathbf{S}$ the set of all $n$-tuples $\mathbf{s}=(s_1,\dots,s_n)\in\{0,1\}^n$ such that
$0\le s_i\le a_i$ for $i=1,\dots, n$, and $(s_i,a_j-s_j)\neq(1,1)$ for every $(i,j)\in Q_B$.
 By convention we assume $a_i=0$ and $s_i=0$ if $i>n$.
Then
\begin{enumerate}\itemsep=0pt
\item[$(i)$] $\tilde{x}[\mathbf{a}]$ can be written as
\begin{gather}\label{df:x[a]}
\left(\prod_{i=1}^n  x_i^{-a_i}\right)\sum_{\mathbf{s}\in \mathbf{S}}\prod_{i=1}^m
x_i^{\sum\limits_{j=1}^n (a_j-s_j)[b_{ij}]_+ + s_j[-b_{ij}]_+}.
\end{gather}

\item[$(ii)$] $\tilde{x}[\mathbf{a}]$ is in the upper cluster algebra $\mathcal{U}$.
\end{enumerate}
\end{Lemma}

\begin{proof}
(i) By Def\/inition \ref{definition:GCC2},
\begin{gather*}
\tilde{x}[\mathbf{a}] =\left(\prod_{i=1}^n  x_i^{-a_i}\right)
\sum_{\mathbf{s}\in \mathbf{S}}\left(\prod_{(i,j)\in Q_{\tilde B}}  x_{i}^{(a_j-s_j)b_{ij}} x_{j}^{s_i(-b_{ji})}\right)\\
 \hphantom{\tilde{x}[\mathbf{a}]}{} =\left(\prod_{i=1}^n  x_i^{-a_i}\right)
\sum_{\mathbf{s}\in \mathbf{S}}\left(\prod_{i,j=1}^m  x_{i}^{(a_j-s_j)[b_{ij}]_+} x_{j}^{s_i[-b_{ji}]_+}\right)\\
\hphantom{\tilde{x}[\mathbf{a}]}{}
=\left(\prod_{i=1}^n  x_i^{-a_i}\right)
\sum_{\mathbf{s}\in \mathbf{S}}\left(\prod_{i,j=1}^m  x_{i}^{(a_j-s_j)[b_{ij}]_+}\right)
\left(\prod_{i,j=1}^m  x_j^{s_i[-b_{ji}]_+}\right)
\\
\hphantom{\tilde{x}[\mathbf{a}]}{}
=\left(\prod_{i=1}^n  x_i^{-a_i}\right)
\sum_{\mathbf{s}\in \mathbf{S}}\left(\prod_{i,j=1}^m  x_{i}^{(a_j-s_j)[b_{ij}]_+}\right)
\left(\prod_{i,j=1}^m  x_i^{s_j[-b_{ij}]_+}\right)
\\
\hphantom{\tilde{x}[\mathbf{a}]}{}
=\left(\prod_{i=1}^n  x_i^{-a_i}\right)
\sum_{\mathbf{s}\in \mathbf{S}}\left(\prod_{i,j=1}^m  x_{i}^{(a_j-s_j)[b_{ij}]_++s_j[-b_{ij}]_+}\right)=\eqref{df:x[a]}.
\end{gather*}

(ii) We introduce $n$ extra variables $x_{m+1},\dots, x_{m+n}$. Let
\begin{gather*}
\mathbb{ZP'}=
\mathbb{ZP}\big[x_{m+1}^{\pm1},\dots,x_{m+n}^{\pm1}\big]=
\mathbb{Z}\big[x_{n+1}^{\pm1},\dots,x_{m+n}^{\pm1}\big]
\end{gather*}
 be the ring of Laurent polynomials in the variables $x_{n+1},\dots,x_{m+n}$. Let
\begin{gather*}
\tilde{B}'=\begin{bmatrix}\tilde{B}\\I_n\end{bmatrix}
 \end{gather*}
 be the $(m+n)\times n$ matrix that encodes a new cluster algebra $\mathcal{A}'$. Assume $a_{n+1}=\cdots=a_{m+n}=0$ and
def\/ine $Q_{\tilde B}'$, $\mathbf{S}'$, $\tilde{x}'[\mathbf{a}]$ for $\mathbb{ZP}'$ similarly as the def\/inition of $Q_{\tilde B}$, $\mathbf{S}$, $\tilde{x}[\mathbf{a}]$ for $\mathbb{ZP}$. (Of course $\mathbf{S}'=\mathbf{S}$, but we use dif\/ferent notation to emphasize that they are for dif\/ferent cluster algebras). So
\begin{gather*}
\tilde{x}'[\mathbf{a}]=\left(\prod_{i=1}^n  x_i^{-a_i}\right)\sum_{{\bf s}\in {\bf S}'}P_{\bf s}, \qquad P_{\bf s}:=\prod_{i=1}^{m+n}
x_i^{\sum\limits_{j=1}^n (a_j-s_j)[b_{ij}]_+ + s_j[-b_{ij}]_+}.
\end{gather*}
We will show that $\tilde{x}'[{\bf a}]$ is in the upper bound
\begin{gather*}
\mathcal{U}'_x=\mathbb{ZP'}\big[{\bf x}^{\pm1}\big]\cap\mathbb{ZP'}\big[{\bf x}_1^{\pm1}\big]\cap\cdots\cap\mathbb{ZP'}\big[{\bf x}_n^{\pm1}\big]
\end{gather*}
(where ${\bf x}=\{x_1,\dots,x_n\}$ and for $1\le k\le n$, the adjacent cluster ${\bf x}_k$ is def\/ined by ${\bf x}_k={\bf x}-\{x_k\}\cup \{x_k'\}$), therefore~$\tilde{x}'[{\bf a}]$ is in the upper cluster algebra~$\mathcal{U}'$, thanks to the fact that $\mathcal{U}'_x=\mathcal{U}'$ when $\tilde{B}'$ is of full rank \cite[Corollary~1.7 and Proposition~1.8]{BFZ}. Then we substitute $x_{n+1}=\cdots=x_{m+n}=1$ and conclude that $\tilde{x}[\mathbf{a}]$ is in the upper cluster algebra $\mathcal{U}$.

Since $\tilde{x}'[\mathbf{a}]$ is obviously in $\mathbb{ZP'}[{\bf x}^{\pm1}]$ from its def\/inition, we only need to show that $\tilde{x}'[\mathbf{a}]$ is in~$\mathbb{ZP'}[{\bf x}_k^{\pm1}]$ for $1\le k\le n$. Again from its def\/inition we see that $\tilde{x}'[\mathbf{a}]$ is in $\mathbb{ZP'}[{\bf x}_k^{\pm1}]$ when $a_k=0$. So we may assume $a_k=1$. Let $N\subset{\bf S}'$ contain those ${\bf s}$ such that $P_{\bf s}$ is not divisible by $x_k$; equivalently,
\begin{gather*}
N=\{{\bf s}\in {\bf S}' \,| \,s_j=a_j\textrm{ if }b_{kj}>0; \; s_j=0 \; \textrm{if} \; b_{kj}<0 \big\}.
\end{gather*}
Then it suf\/f\/ices to show that $\sum\limits_{{\bf s}\in N}P_{\bf s}$ is divisible by $A$, where
\begin{gather*}
A=x_k'x_k=\prod_{i=1}^{m+n}x_i^{[b_{ik}]_+}+\prod_{i=1}^{m+n}x_i^{[-b_{ik}]_+}.
\end{gather*}
Write $N$ into a partition $N=N_0\cup N_1$ where $N_0=\{{\bf s}\in N | s_k=0\}$, $N_1=\{{\bf s}\in N | s_k=1\}$.
Def\/ine $\varphi\colon N_0\to N_1$ by
\begin{gather*}
\varphi(s_1,\dots,s_{k-1},0,s_{k+1},\dots,s_n)=(s_1,\dots,s_{k-1},1,s_{k+1},\dots,s_n).
\end{gather*}
Then $\varphi$ is a well-def\/ined bijection because in the def\/inition of $N$ there is no condition imposed on~$s_k$.
Thus
\begin{gather*}
\sum_{{\bf s}\in N}P_{\bf s} =\sum_{{\bf s}\in N_0}P_{\bf s}+\sum_{{\bf s}\in N_1}P_{\bf s}=
\sum_{{\bf s}\in N_0}(P_{\bf s}+P_{\varphi({\bf s})})\\
\hphantom{\sum_{{\bf s}\in N}P_{\bf s}}{}
=\prod_{i=1}^{m+n}\prod_{\stackrel{1\le j\le n}{j\neq k}}
x_i^{(a_j-s_j)[b_{ij}]_+ + s_j[-b_{ij}]_+}
\left(\prod_{i=1}^{m+n}
x_i^{[b_{ik}]_+}+\prod_{i=1}^{m+n}
x_i^{[-b_{ik}]_+}\right)\\
\hphantom{\sum_{{\bf s}\in N}P_{\bf s}}{}
 =\prod_{i=1}^{m+n}\prod_{\stackrel{1\le j\le n}{j\neq k}}
x_i^{(a_j-s_j)[b_{ij}]_+ + s_j[-b_{ij}]_+}
A
\end{gather*}
is divisible by $A$.
\end{proof}

Below is the main theorem of the section.

\begin{proof}[Proof of Theorem~\ref{thm:x[a]}]
It follows immediately from Lemma~\ref{multiplicative}(i) and Lemma~\ref{lemma:01}(ii).
Indeed, we use Lemma \ref{multiplicative}(i) to factor $\tilde{x}[{\bf a}]$ (for ${\bf a}\in\mathbb{Z}_{\ge0}^n$) into the product of a usual monomial (which is in $\mathcal{U}$) and those $\tilde{x}[{\bf a}']$'s where all entries of ${\bf a}'$ are $0$ or $1$ (which is also in~$\mathcal{U}$ by  Lemma~\ref{lemma:01}(ii).
\end{proof}

\begin{Remark}
We claim that if $\Sigma$ is acyclic (i.e., $Q_B$ is acyclic), then $\big(\prod\limits_{i=1}^n  x_i^{a_i}\big)\tilde{x}[\mathbf{a}]$ is not divisible by~$x_k$ for any $1\le k\le n$, i.e.,  there exists  $\mathbf{s}=(\mathbf{s}_1,\dots,\mathbf{s}_n)\in S_{\rm gcc}$  such that
\begin{gather*}
\prod_{(i,j)\in Q_{\tilde B}}  x_{i}^{b_{ij}|\bar{\mathbf{s}}_j|} x_{j}^{-b_{ji}|\mathbf{s}_i|}
\end{gather*}
is not divisible by~$x_k$.
Fix $k$ such that $1\le k\le n$. We need to f\/ind ${\bf s}$ such that
\begin{gather*}
 |\bar{\mathbf{s}}_j|=0 \quad \text{if} \ \ (k,j)\in Q_{\tilde B} \qquad \text{and} \qquad |\mathbf{s}_i|=0 \quad \text{if} \ \ (i,k)\in Q_{\tilde B}.
\end{gather*}
This condition is equivalent to
\begin{gather*}
|\bar{\mathbf{s}}_j|=0 \quad \text{if} \ \ (k,j)\in Q_B \qquad \text{and} \qquad |\mathbf{s}_i|=0 \quad \text{if} \quad (i,k)\in Q_B,
\end{gather*}
because ${\bf s}_i=()$ for $i>n$ by convention.
Such  an ${\bf s}$ can be constructed as follows: since  the initial seed is acyclic, $Q_B$ has no oriented cycles, therefore it determines a partial order where $i\prec j$ if there is a (directed) path from~$i$ to~$j$ in~$Q_B$. Def\/ine
$\mathbf{s}_l\in\{0,1\}^{[a_l]_+}$ for $l=1,\dots,n$ as
\begin{gather*}
\mathbf{s}_l=\begin{cases}
(1,\dots,1),& \textrm{if} \ \ k\prec l,\\
(0,\dots,0),& \textrm{otherwise}.
\end{cases}
\end{gather*}
Note that the construction of ${\bf s}$ depends on $k$. To check that ${\bf s}$ is in $S_{\rm gcc}$, we need to show that $\mathbf{s}_i\cdot\bar{\mathbf{s}}_j=0$ for every $(i,j)\in Q_B$. This can be proved by contradiction as follows.
Assume $\mathbf{s}_i\cdot\bar{\mathbf{s}}_j\neq 0$ for some $(i,j)\in Q_B$.
Then there exists $r\le\min([a_i]_+,[a_j]_+)$ such that $s_{i,r}=1$, $s_{j,r}=0$. By our choice of $\mathbf{s}$, we have $k\prec i$, $k\not\prec j$, therefore $i\not\prec j$. This contradicts with the assumption that $(i,j)\in Q_B$.
\end{Remark}

\section{Non-acyclic rank 3 cluster algebras}\label{section4}
In this section, we consider non-acyclic skew-symmetric rank 3 cluster algebras of geometric type.

\subsection[Def\/inition and properties of $\tau$]{Def\/inition and properties of $\boldsymbol{\tau}$}\label{section4.1}

Let  $m\geq 3$ be an integer, $\tilde{B}=(b_{ij})$ be an $m \times 3$ matrix whose principal part $B$ is skew-symmetric and non-acyclic (i.e., $b_{12}$, $b_{23}$, $b_{31}$ are of the same sign).

Def\/ine the map
\begin{gather*}
\tau(B):=(|b_{23}|,|b_{31}|,|b_{12}|)\in \mathbb{Z}^3_{\geq 0}.
\end{gather*}

\begin{Definition}
 Let $(x,y,z) \in \mathbb{R}^3$. We def\/ine a partial ordering ``$\le$'' on  $\mathbb{R}^3$ by $(x,y,z) \leq (x',y',z')$ if and only if $x \leq x'$, $y \leq y'$, and $z \leq z'$.
\end{Definition}

Def\/ine three involutary functions $\mu_1,\mu_2,\mu_3\colon \mathbb{R}^3\to\mathbb{R}^3$ as follows
\begin{gather*}
\mu_1\colon \ (x,y,z)\mapsto (x,xz-y,z), \qquad \mu_2\colon \ (x,y,z) \mapsto  (x,y,xy-z),
\\
\mu_3\colon \  (x,y,z) \mapsto(yz-x,y,z).
\end{gather*}
Let $\Gamma$ be the group generated by $\mu_1$, $\mu_2$, and $\mu_3$.  It follows that the $\Gamma$-orbit of $\tau(B)$ is identical to the set $\{\tau(B')$: $B'$ is in the mutation class of~$B \}$.

 Consider the following two situations.
\begin{enumerate}\itemsep=0pt
\item[(M1)] $(a,b,c)\leq \mu_i(a,b,c)$ for all $i=1,2,3$.
\item[(M2)] $(a,b,c)\leq \mu_i(a,b,c)$ for precisely two indices~$i$.
\end{enumerate}

The following statement follows from \cite[Theorem~1.2, Lemma~2.1]{BBH} for a coef\/f\/icient free cluster algebra. The result holds more generally for a cluster algebra of geometric type since the coef\/f\/icients do not af\/fect whether or not the seed is acyclic.
\begin{Theorem}
\label{cyclic}
Let $B$ be as above, $(a,b,c)=\tau(B)$. Then the following are equivalent:
\begin{itemize}\itemsep=0pt
\item[\rm(i)]$\mathcal{A}({\bf x},{\bf y},B)$ is non-acyclic.

\item[\rm(ii)] $a,b,c \geq 2$ and $abc+4\ge a^2+b^2+c^2$.

\item[\rm(iii)] $a,b,c\ge 2$ and there exists a unique triple in the $\Gamma$-orbit of $\tau(B)$ that satisfies~{\rm(M1)}.
\end{itemize}
Moreover, under these conditions, triples in the $\Gamma$-orbit of $\tau(B)$ not satisfying~{\rm(M1)} must sa\-tis\-fy~{\rm(M2)}.
\end{Theorem}

Assume $\mathcal{A}({\bf x},{\bf y},B)$ is non-acyclic. We call the unique triple in Theorem \ref{cyclic}(iii) the \emph{root} of the $\Gamma$-orbit of~$\tau(B)$.

\begin{Lemma}\label{lemma:chain inequality}
Let $B$ be as above, $(a,b,c)=\tau(B)$.
\begin{enumerate}\itemsep=0pt
\item[$(i)$] The root $(a,b,c)$ of the $\Gamma$-orbit of $\tau(B)$ is the minimum of the $\Gamma$-orbit, i.e., $(a,b,c)\le(a',b',c')$ for any $(a',b',c')$ in the $\Gamma$-orbit.

\item[$(ii)$] If $\mu_k(a,b,c)=(a,b,c)$ for any $k$, then $a=b=c=2$ and $(2,2,2)$ is the unique triple in the $\Gamma-$orbit of~$\tau(B)$. Therefore $(2,2,2)$ must also be the root.

\item[$(iii)$] Assume that $\tau(B)$ is the root of the $\Gamma$-orbit of $\tau(B)$,  $i_1,\dots,i_l\in\{1,2,3\}$, $i_s\neq i_{s+1}$ for $1\le s\le l-1$.
Let $B_0=B$, $B_j=\mu_{i_j}(B_{j-1})$ for $j=1,\dots,l$. Then
\begin{gather*}\tau(B_0)\le\tau(B_1)\le\cdots\le\tau(B_l).
\end{gather*}
\end{enumerate}
\end{Lemma}

\begin{proof}
(i) If $\tau(B')$ is the minimal triple in the $\Gamma$-orbit of $\tau(B)$ it satisf\/ies~(M1), so by Theo\-rem~\ref{cyclic}(iii) the root is the unique triple satisfying~(M1) and $B'$ must be the root.
For~(ii): $a=b=c=2$ follows easily from Theorem~\ref{cyclic}(ii), and the uniqueness claim follows from the def\/inition of~$\mu_k$.
 (iii)~is obvious if  $\tau(B)=(2,2,2)$. If not, assume~(iii) is false, then there are $1\le j<j'\le l$ such that
 \begin{gather*}
 \tau(B_j)<\tau(B_{j+1})=\tau(B_{j+2})=\cdots=\tau(B_{j'})>\tau(B_{j'+1}).
 \end{gather*}
By (ii) it must be that $j+1=j'$ so that  $\tau(B_j)< \tau(B_{j+1}) > \tau(B_{j+2})$, but $\tau(B_{j+1})$ does not satisfy either~(M1) or~(M2). This is impossible by Theorem~\ref{cyclic}(iii).
\end{proof}

\subsection[Grading on~$\mathcal{A}$]{Grading on $\boldsymbol{\mathcal{A}}$}\label{section4.2}

We now adapt the grading introduced in \cite[Def\/inition 3.1]{G} to the geometric type.
\begin{Definition}
A graded seed is a quadruple $({\bf x},{\bf y},B,G)$ such that
\begin{enumerate}\itemsep=0pt
\item[(i)] $({\bf x},{\bf y},B)$ is a seed of rank $n$, and
\item[(ii)] $G=[g_1,\dots,g_n]^T\in \mathbb{Z}^n$ is an integer column vector such that $BG=0$.
\end{enumerate}
\end{Definition}

Set $\deg_G(x_i)=g_i$ and $\deg\big(x_i^{-1}\big)=-g_i$ for $i\le n$, and set $\deg(y_j)=0$ for all $j$. Extend the grading additively to Laurent monomials hence to the cluster algebra $\mathcal{A}({\bf x},{\bf y},B)$. In~\cite{G} it is proved that under this grading every exchange relation is homogeneous and thus the grading is compatible with mutation.

\begin{Theorem}[\protect{\cite[Corollary 3.4]{G}}]
The cluster algebra $\mathcal{A}({\bf x},{\bf y},B)$ under the above grading is a~$\mathbb{Z}$-graded algebra.
\end{Theorem}

The following two propositions come from the work in \cite{S1} to show that rank three non-acyclic cluster algebras have no maximal green sequences.
\begin{Theorem}[\protect{\cite[Proposition 2.2]{S1}}]
Suppose that $B$ is a $3\times 3$ skew-symmetric non-acyclic matrix. Then the $($column$)$ vector $G=\tau(B)^T$ satisfies $BG=0.$
\end{Theorem}

\begin{Lemma}\label{lemma: grading}
For any graded seed $({\bf x}',B',G')$ in the mutation class of our initial graded seed, we have $G'= \tau(B')^T$.
\end{Lemma}

\begin{proof}
This result follows from \cite[Lemma~2.3]{S1}, but its proof is short, so we reproduce it here. We use induction on mutations. Suppose that $G'= \tau(B')^T$ for a given graded seed $({\bf x}',B',G')$. Let $({\bf x}'',B'',G'')$ be the graded seed obtained by taking mutation $\mu_1$ from $({\bf x}',B',G')$.  Then $x_1''=((x_2')^{|b'_{12}|}+(x_3')^{|b'_{13}|})/x_1'$, so its degree is $g_2' |b'_{12}|-g_1'$. By induction we have
\begin{gather*}
g_2' |b'_{12}|-g_1' =|b'_{31}| |b'_{12}|-|b'_{23}| = |b'_{31}b'_{12}-b'_{23}| =|b''_{23}|,
\end{gather*}
so we get $G''= \tau(B'')^T$. The cases of~$\mu_2$ and~$\mu_3$ are similar.
\end{proof}

In the rest of Section~\ref{section4} we assume that $B$ is a $3\times 3$ skew-symmetric non-acyclic matrix and $G=\tau(B)^T$. In other words, if $\tau(B)=(b,c,a)$ then $\deg(x_1)=b$, $\deg(x_2)=c$, $\deg(x_3)=a$ for any seed $({\bf x}, {\bf y}, B)$. If we look at the quiver associated to our exchange matrix we see that the degree of the cluster variable $x_i$ is the number of arrows between the other two mutable vertices in~$Q_B$. Furthermore, this is a~canonical grading for non-acyclic rank three cluster algebras since regardless of your choice of initial seed the grading imposed on~$\mathcal{A}$ is the same.

\subsection[Construction of an element in $\mathcal{U}\setminus\mathcal{A}$]{Construction of an element in $\boldsymbol{\mathcal{U}\setminus\mathcal{A}}$}\label{section4.3}
Here we shall prove the following main theorem of the paper.

We give $\mathcal{A}$ the $\mathbb{Z}$-grading as in Section~\ref{section4.2}. By Theorem~\ref{cyclic}, we can assume that the initial seed
\begin{gather*}
\Sigma=(\{x_1,x_2,x_3\},{\bf y},B,G),\qquad\textrm{where} \qquad B= \left[ \begin{matrix}
0 & a &-c \\
-a & 0 & b \\
c& -b & 0
\end{matrix}  \right], \qquad G=\left[ \begin{matrix}
b \\ c \\ a  \end{matrix} \right]
\end{gather*}
satisf\/ies $a,b,c\ge2$ and (M1). Then $\deg x_{1}=b$, $\deg x_{2}=c$, $\deg x_{3}=a$, and we let $x_{4},\dots,x_{m}$ be of degree~0.  Furthermore, by permuting the indices if necessary, we assume
$a\ge b\ge c$.

Def\/ine six degree-0 elements as follows
\begin{gather*}
\alpha_i^\pm := \prod_{j=4}^m x_j^{[\pm b_{ji}]_+}, \qquad i=1,2,3.
\end{gather*}
Looking at the seeds neighboring $\Sigma$ we obtain three new cluster variables. Namely,
\begin{gather*}
z_1=\frac{\alpha_1^-x_2^a+\alpha_1^+x_3^c}{x_1},
\qquad z_2=\frac{\alpha_2^+x_1^a+\alpha_2^-x_3^b}{x_2}, \qquad z_3=\frac{\alpha_3^-x_1^c+\alpha_3^+x_2^b}{x_3},
 \end{gather*} which have degree $ac-b$, $ab-c$, and $bc-a$, respectively.

We also use the following theorem, whose proof in~\cite{MM} is for coef\/f\/icient free cluster algebras but clearly applies to cluster algebras of geometric type.
\begin{Theorem}[\protect{\cite[Proposition 6.1.2]{MM}}]
\label{coprime}
If $\mathcal{A}$ is a rank three skew-symmetric non-acyclic cluster algebra, then~$\mathcal{A}$ is totally coprime.
\end{Theorem}

Now we consider a special element in $\mathbb{Q}(x_1,x_2,x_3)$:
\begin{gather*}
Y:=\tilde{x}[(1,0,1)]/x_2^b =\frac{\alpha_1^-\alpha_3^-x_1^cx_2^{a-b}+\alpha_1^-\alpha_3^+x_2^a+\alpha_1^+\alpha_3^+x_3^c}{x_1x_3}.
\end{gather*}
We shall prove Theorem~\ref{main-cor} by showing that the element constructed above is in $\mathcal{U}\setminus\mathcal{A}$.

\begin{proof}[Proof of Theorem \ref{main-cor}]
We f\/irst show that $Y\in \mathcal{U}$.
Note that $\mathcal{A}$ is a rank 3 cluster algebra so it is totally coprime by Theorem~\ref{coprime}. It then suf\/f\/ices to show that $Y\in \mathcal{U}_{\bf x}$ by Theorem~\ref{UB}. Clearly $Y\in\mathbb{ZP}[x_1^{\pm 1},x_2^{\pm1},x_3^{\pm 1}]$, where
$\mathbb{ZP}=\mathbb{Z}[x_4^{\pm1},\dots,x_m^{\pm1}]$. Also,
\begin{gather*}
Y=\frac{\alpha_3^+z_1^c+\alpha_1^-\alpha_3^-\big(\alpha_1^-x_2^a+\alpha_1^+x_3^c\big)^{c-1} x_2^{a-b}}{z_1^{c-1}x_3} \in \mathbb{ZP}\big[z_1^{\pm 1},x_2^{\pm 1},x_3^{\pm 1}\big],\\ Y=\frac{\alpha_1^-\alpha_3^-x_1^c\!\big(\alpha_2^+x_1^a\!+\!\alpha_2^-x_3^b\big)^{a{-}b}\!z_2^b\!+\!\alpha_1^+\alpha_3^+z_2^ax_3^c\!+\!\alpha_1^-\alpha_3^+
\!\big(\alpha_2^+x_1^a\!+\!\alpha_2^-x_3^b\big)^{a}}{x_1z_2^ax_3}
 \in  \mathbb{ZP}\big[x_1^{\pm 1}\!,z_2^{\pm 1}\!,x_3^{\pm 1}\big],\\
  Y=\frac{\alpha_1^-x_2^{a-b}z_3^c+\alpha_1^+\alpha_3^+\big(\alpha_3^-x_1^c+\alpha_3^+x_2^b\big)^{c-1}}{x_1z_3^{c-1}} \in\mathbb{ZP}\big[x_1^{\pm 1},x_2^{\pm 1},z_3^{\pm 1}\big].
\end{gather*}
Therefore we conclude that $Y\in \mathcal{U}_{\bf x}$.

Next, we show that $Y\notin \mathcal{A}$. With respect to our grading of $\mathcal{A}$, $Y$ is homogeneous of degree $ac-b-a.$

Combining Lemmas \ref{lemma:chain inequality} and~\ref{lemma: grading} we have already shown that the degree of cluster variables is non-decreasing as we mutate away from our initial seed. We use this fact to explicitly prove that all cluster variables, except possibly $x_1$, $x_2$, $x_3$, $z_3$, have degree strictly larger than $ac-b-a$.  Indeed, let $x$ be a cluster variable such that $x\neq x_1,x_2,x_3,z_3$. Then $x$ can be written as
\begin{gather*}
x=\mu_{i_l}\cdots\mu_{i_2}\mu_{i_1}(x_k),\quad i_1,\dots,i_l,k\in\{1,2,3\}, \qquad \text{and}\qquad i_s\neq i_{s+1} \quad \text{for} \ \ 1\le s\le l-1.
\end{gather*}
Let $B_0=B$, $B_j=\mu_{i_j}(B_{j-1})$ for $j=1,\dots,l$. Let $\tau(B_j)=(b_j,c_j,a_j)$ for $j=0,\dots,l$. By Lemma~\ref{lemma:chain inequality}, we have that
\begin{gather*}
(b,c,a)=(b_0,c_0,a_0)\le(b_1,c_1,a_1)\le\cdots\le(b_l,c_l,a_l).
\end{gather*}
We prove the claim in the following f\/ive cases.

\textit{Case $k=1$}. We let $r>0$ be the smallest integer such that $i_r=k$ (which exists since $x\neq x_1,x_2,x_3$). It suf\/f\/ices to prove that the degree of $w=\mu_{i_r}\cdots\mu_{i_2}\mu_{i_1}(x_k)$ is larger than $ac-b-a$. Indeed,
\begin{gather*}
\deg w=b_r=a_{r-1}c_{r-1}-b_{r-1}=a_{r-1}c_{r-1}-b\ge ac-b>ac-b-a.
\end{gather*}

\textit{Case $k=2$}. The above proof still works:
\begin{gather*}
\deg w=c_r=a_{r-1}b_{r-1}-c_{r-1}=a_{r-1}b_{r-1}-c\ge ab-c>ac-b-a.
\end{gather*}

\textit{Case  $k=3$, $i_1=1$}. We let $r>0$ be the smallest integer such that $i_r=k$. Then $(b_1,c_1,a_1)=(ac-b,c,a)$, and
\begin{gather*}
\deg w=a_r=b_{r-1}c_{r-1}-a_{r-1}=b_{r-1}c_{r-1}-a\ge b_1c_1-a=(ac-b)c-a>ac-b-a.
\end{gather*}

\textit{Case  $k=3$, $i_1=2$}. Similar to the above case, $(b_1,c_1,a_1)=(b,ab-c,a)$,
\begin{gather*}
\deg w=a_r\ge b_1c_1-a=b(ab-c)-a\ge b^2(a-1)-a>c(a-1)-a\ge ac-b-a.
\end{gather*}
 The second and fourth inequality follow from our assumption that $b \geq c$.

\textit{Case  $k=3$, $i_1=3$}. Then $(b_1,c_1,a_1)=(b,c,bc-a)$, $i_2\neq 3$. We let $r\ge3$ be the smallest integer such that $i_r=k$ (which exists since $x\neq z_3$).  It suf\/f\/ices to prove that the degree of $w=\mu_{i_r}\cdots\mu_{i_2}\mu_{i_1}(x_k)$ is larger than $ac-b-a$.

If $i_2=1$, then $(b_2,c_2,a_2)=(c(bc-a)-b,c,bc-a)$,
\begin{gather*}
\begin{split}
&\deg w =a_r=b_{r-1}c_{r-1}-(bc-a)\ge b_2c_2-(bc-a)=(c(bc-a)-b)c-(bc-a)
\\
& \hphantom{\deg w }{} =(c^2-2)(bc-a)-a\ge(c^2-2)a-a=a(c^2-3)\ge a(c-1)>ac-b-a.
\end{split}
\end{gather*}

If $i_2=2$, then $(b_2,c_2,a_2)=(b,(bc-a)b-c,bc-a)$,
\begin{gather*}
\deg w\ge b((bc-a)b-c)-(bc-a)\ge (c(bc-a)-b)c-(bc-a)>ac-b-a.
\end{gather*}
The second inequality above follows from our assumption that $b \geq c$. One application gives the inequality $(bc-a)b-c\ge c(bc-a)-b$ and a second application replaces a factor of $b$ in the f\/irst term with a factor of~$c$. The last inequality follows from our computation in the previous case.

This completes the proof of the claim that all cluster variables, except possibly $x_1$, $x_2$,~$x_3$,~$z_3$, have degree strictly larger than $ac-b-a$.

Therefore it is suf\/f\/icient to check that $Y$ cannot be written as a linear combination of products of the cluster variables $x_1$, $x_2$, $x_3$, and~$z_3$ since all other cluster variable have larger degree than~$Y$. This is clear as the numerator of~$Y$ is irreducible and none of these cluster variables have a~factor of~$x_1$ in the denominator.
\end{proof}

\section{Dyck path formula}\label{section5}

 In this section, we give the Dyck path construction of $\tilde{x}[{\bf a}]$ that is equivalent to Section~\ref{section3} (under a mild condition that there is no isolated vertex in the digraph $Q_{\tilde B}$). It is in fact our original def\/inition of $\tilde{x}[{\bf a}]$, and appears naturally in the attempt of  generalizing greedy bases for cluster algebras of rank~2 in~\cite{LLZ} and the construction of bases of type~$A$ cluster algebras in~\cite{B9} to more general cases.

Let $(a_1, a_2)$ be a pair of nonnegative integers. Let $c=\min(a_1,a_2)$.
The \emph{maximal Dyck path} of type $a_1\times a_2$, denoted by $\mathcal{D}=\mathcal{D}^{a_1\times a_2}$, is a lattice path
from $(0, 0)$  to $(a_1,a_2)$ that is as close as possible to the diagonal joining $(0,0)$ and $(a_1,a_2)$, but never goes above it. A~\emph{corner} is a~subpath consisting of a~horizontal edge followed by a~vertical edge.

\begin{Definition}\label{corner-first}
Let $\mathcal{D}_1$ (resp.~$\mathcal{D}_2$) be the set of horizontal (resp.~vertical) edges of a maximal Dyck path $\mathcal{D}=\mathcal{D}^{a_1\times a_2}$. We label $\mathcal{D}$ with the \emph{corner-first index} in the following sense:
\begin{enumerate}\itemsep=0pt
\item[(a)] edges in $\mathcal{D}_1$ are indexed as $u_1,\dots,u_{a_1}$ such that $u_i$ is the horizontal edge of the $i$-th corner for $i\in [1,c]$ and $u_{c+i}$ is the $i$-th of the remaining horizontal ones for $i\in [1, a_1-c]$,

\item[(b)] edges in $\mathcal{D}_2$ are indexed as $v_1,\dots, v_{a_2}$ such that $v_i$ is the vertical edge of the $i$-th corner for $i\in [1,c]$ and $v_{c+i}$ is the $i$-th  of the remaining vertical ones for $i\in [1, a_2-c]$.
\end{enumerate}
 Here we count corners from bottom left to top right, count vertical edges from bottom to top, and count horizontal edges from left to right.)
\end{Definition}

\begin{Definition}\label{local_compatibility}
Let $S_1\subseteq \mathcal{D}_1$ and  $S_2 \subseteq \mathcal{D}_2$. We say that $S_1$ and $S_2$ are \emph{locally compatible} with respect to $\mathcal{D}$ if and only if no horizontal edge in $S_1$ is the immediate predecessor of any vertical edge in $S_2$ on $\mathcal{D}$. In other words, the subpath $S_1\cup S_2$ contains no corners.
\end{Definition}

\begin{Remark}
Notice that we have: $|\mathcal{P(D}_1)\times\mathcal{P(D}_2)|$ possible pairs for $(S_1,S_2)$, where $\mathcal{P(D}_1)$ denotes the power set of $\mathcal{D}_1$ and $\mathcal{P(D}_2)$ denotes the power set of $\mathcal{D}_2$. Further, when either $S_1$ or $S_2 = \varnothing$, any arbitrary choice of the other will yield local compatibility by our above def\/inition.
\end{Remark}

\begin{figure}[h!]
\begin{enumerate}[(i)]
\itemsep=0pt
\item
For each corner,
\includegraphics{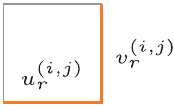},
\includegraphics{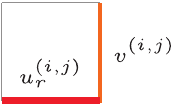},
or
\includegraphics{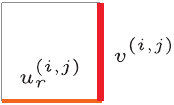}
are allowed,
but \includegraphics{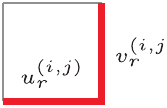} is not.

\item
\includegraphics{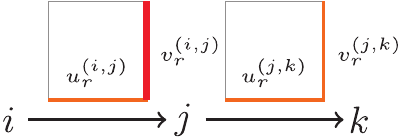} \ or \ \includegraphics{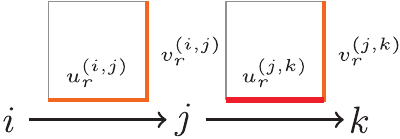}

\item
\includegraphics{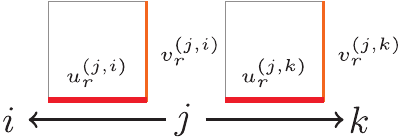} \ or \ \includegraphics{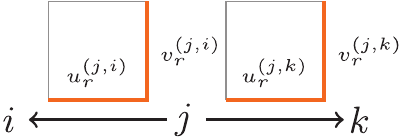}

\item
\includegraphics{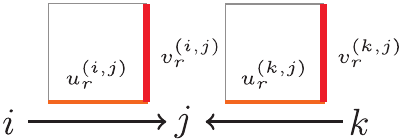} \ or \ \includegraphics{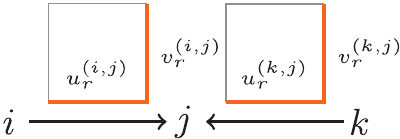}
\end{enumerate}
\caption{Pictorial descriptions for Def\/inition~\ref{definition:GCC}(i)--(iv).}
\label{GCCillust}
\end{figure}

\begin{Definition}\label{definition:GCC}
Let $\mathbf{a}=(a_i)\in \mathbb{Z}^n$. By convention we assume $a_i=0$ for $i>n$.
For each pair $(i,j)\in Q_{\tilde B}$ (def\/ined at the beginning of Section~\ref{section3}), denote $\mathcal{D}^{[a_i]_+\times [a_j]_+}$ by $\mathcal{D}^{(i,j)}$ . We label  $\mathcal{D}^{(i,j)}$ with the corner-f\/irst index (Def\/inition~\ref{corner-first}), whose horizontal edges are denoted $u_{1}^{(i,j)},\dots,u_{[a_{i}]_{+}}^{(i,j)}$ and vertical edges are denoted by $v_{1}^{(i,j)},\dots,v_{[a_{j}]_{+}}^{(i,j)}$.
We say that the collection
\begin{gather*}
\big\{S_\ell^{(i,j)} \subseteq\mathcal{D}^{(i,j)}_\ell\,\big|\,  (i,j)\in Q_{\tilde B},\; \ell\in\{1,2\}\big\}
\end{gather*}
 is \emph{globally compatible} if and only if
\begin{enumerate}[(i)]\itemsep=0pt

 \item $S^{(i,j)}_1$ and $S^{(i,j)}_2$ are \emph{locally compatible} with respect to $\mathcal{D}^{(i,j)}$ for all $(i,j)\in Q_{\tilde B}$;

\item if $(i,j)$ and $(j,k)$ are in $Q_{\tilde B}$, then $v^{(i,j)}_r\in S^{(i,j)}_2\Longleftrightarrow u^{(j,k)}_r\not\in S^{(j,k)}_1$ for all $r\in[1,[a_{j}]_{+}]$;

\item if $(j,k)$ and $(j,i)$ are in $Q_{\tilde B}$, then $u^{(j,k)}_r \in S^{(j,k)}_1\Longleftrightarrow u^{(j,i)}_r \in S^{(j,i)}_1$ for all $r\in[1,[a_{j}]_{+}]$;

\item if $(i,j)$ and $(k,j)$ are in $Q_{\tilde B}$, then $v^{(i,j)}_r \in S^{(i,j)}_2\Longleftrightarrow v^{(k,j)}_r \in S^{(k,j)}_2$ for all $r\in[1,[a_{j}]_{+}]$.
\end{enumerate}
We say that the collection is \emph{quasi-compatible} if it satisf\/ies (ii), (iii), and (iv). So a globally compatible collection is also a quasi-compatible collection.
For an illustration of the conditions, see Fig.~\ref{GCCillust}.
\end{Definition}

\begin{Proposition}\label{prop:GCCxa} Same notation as in Definition~{\rm \ref{definition:GCC}}. Assume that there is no isolated vertex in the digraph $Q_{\tilde B}$. Then
\begin{gather*}
\tilde{x}[\mathbf{a}]=\left(\prod_{l=1}^n  x_l^{-a_l}\right)
\sum\left(\prod_{(i,j)\in Q_{\tilde B}}  x_{i}^{b_{ij}\big|S^{(i,j)}_2\big|} x_{j}^{-b_{ji}\big|S^{(i,j)}_1\big|}\right) ,
\end{gather*}
where the sum runs over all globally compatible collections $($abbreviated GCCs$)$. A similar formula holds for~$z[\mathbf{a}]$ if the sum runs over all quasi-compatible collections.
\end{Proposition}

\begin{proof}
We def\/ine a mapping that sends a globally compatible collection
\begin{gather*}
\big\{S_\ell^{(i,j)} \subseteq\mathcal{D}^{(i,j)}_\ell\,\big|\,  (i,j)\in Q_{\tilde B},\; \ell\in\{1,2\}\big\},
\end{gather*}
to a sequence of sequences ${\bf s}=({\bf s}_1,\dots,{\bf s}_n)$ as follows. For $1\le j\le n$ and $1\le r\le [a_j]_+$, the $r$-th number $s_{j,r}$ of the sequence $\mathbf{s}_j=(s_{j,1}, s_{j,2},\dots,s_{j,[a_j]_+})\in\{0,1\}^{[a_j]_+}$ is determined by
\begin{enumerate}[(i)]\itemsep=0pt
\item if $(i,j)\in Q_{\tilde B}$, then $v^{(i,j)}_r\in S^{(i,j)}_2$ if and only if $s_{j,r}=0$;
\item if $(j,k)\in Q_{\tilde B}$, then $u^{(j,k)}_r\in S^{(j,k)}_1$ if and only if $s_{j,r}=1$.
\end{enumerate}
Since we assume that~$Q_{\tilde B}$ has no isolated vertices and the collection $\{S_\ell^{(i,j)}\}$ are globally compatible,  all~$s_{j,r}$ are well-def\/ined. Noting that the number of corners in $\mathcal{D}^{(i,j)}$ is $\min([a_i]_+,[a_j]_+)$ and that we use the corner-f\/irst index, it is easy to conclude that ${\bf s}$ is indeed in~$S_{\rm gcc}$, and that the mapping is bijective and the exponents that appear in $\tilde{x}[{\bf a}]$ work out correctly.

The formula for~$z[a]$ is proved similarly.
\end{proof}

\begin{Remark}
Note that the local compatibility condition  in Def\/inition~\ref{definition:GCC}(i) is weaker than the compatibility condition for greedy elements given in \cite{LLZ}. Even for rank 2 cluster algebras, not all cluster variables are of the form $\tilde{x}[\mathbf{a}]$. As an example,
\begin{gather*}
B=\begin{bmatrix}
0&2\\-2&0
\end{bmatrix}.
\end{gather*}
Then by Lemma \ref{multiplicative}(i),
\begin{gather*}
\tilde{x}[(1,2)]=\tilde{x}[(1,1)]\tilde{x}[(0,1)]=\left(\frac{1+x_1^2+x_2^2}{x_1x_2}\right)
\left(\frac{1+x_1^2}{x_2}\right)=\frac{1+2x_1^2+x_1^4+x_2^2+x_1^2x_2^2}{x_1x_2^2},
\end{gather*}
which is not equal to the following cluster variable (which is a greedy element)
\begin{gather*}
x[(1,2)]=\frac{1+2x_1^2+x_1^4+x_2^2}{x_1x_2^2}.
\end{gather*}
Other $\tilde{x}[{\bf a}]$ are not equal to this cluster variable either, since their denominators do not match. It is illustrating to also compare with the standard monomial basis element
\begin{gather*}
z[(1,2)]=\left(\frac{1+x_2^2}{x_1}\right)\left(\frac{1+x_1^2}{x_2}\right)^2=\frac{1+2x_1^2+x_1^4+x_2^2+2x_1^2x_2^2+x_1^4x_2^2}{x_1x_2^2}.
\end{gather*}
Similar as the standard monomial basis, $\{\tilde{x}[{\bf a}]\}$ also depends on the initial cluster in general.

For interested readers, we explain the above example using pictures:  in Fig.~\ref{figure1},  all 8 possible collections of edges in~$\mathcal{D}^{1\times 2}$ together with their contributions to the numerator of~$z[(1,2)]$ are listed; among them, the f\/irst 6 are GCCs and contribute to the numerator of~$\tilde{x}[(1,2)]$, and the f\/irst 5 are compatible pairs in the def\/inition of greedy elements and contribute to the numerator of~$x[(1,2)]$.
\begin{figure}[h!]\centering
\includegraphics{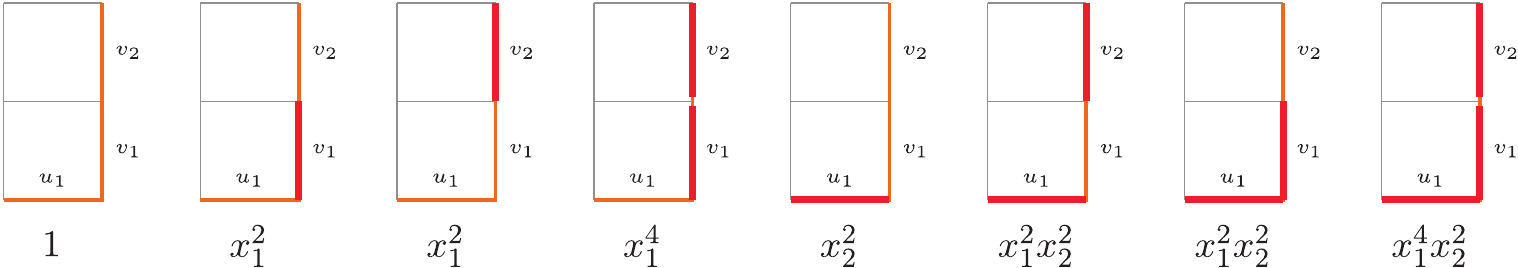}
\caption{Qausi-compatible collections for ${\bf a}=(1,2)$.}
\label{figure1}
\end{figure}
\end{Remark}

 \section{Proof of Theorem \ref{thm:basis}}\label{section6}

 In this section, we give two proofs of the fact that $\{\tilde{x}[\textbf{a}]\}_{\textbf{a}\in\mathbb{Z}^n}$ constructed from an {\it acyclic seed} (i.e., {\it $Q_B$~is acyclic}) form a basis of an acyclic cluster algebra $\mathcal{A}$ of geometric type, using dif\/ferent approaches. The main idea behind both proofs is to compare $\{\tilde{x}[\textbf{a}]\}_{\textbf{a}\in\mathbb{Z}^n}$ with the standard monomial basis $\{z[\textbf{a}]\}_{\textbf{a}\in\mathbb{Z}^n}$ which is known to have the desired property.
 The f\/irst approach is to consider certain orders on Laurent monomials and use the fact that $\tilde{x}[{\bf a}]$ and $z[{\bf a}]$  have the same lowest Laurent monomial to draw the conclusion. The second approach is to use the multiplicative property of $\tilde{x}[{\bf a}]$ and~$z[{\bf a}]$ and the combinatorial descriptions of these elements.

\subsection{The f\/irst proof}\label{section6.1}

 Note that $Q_B$ is acyclic by assumption.
If $Q_B$ has only one vertex, then $\tilde{x}[{\bf a}]=z[{\bf a}]$ and we are done.
If~$Q_B$ is not weakly connected (recall that a digraph is weakly connected if replacing all of its directed edges with undirected edges produces a connected undirected graph), we can factorize $\tilde{x}[{\bf a}]$ (resp.\ $z[{\bf a}]$) into a product of those with smaller ranks, each corresponds to a~connected component (see Lemma \ref{multiplicative}(ii)). Thus we can reduce to the situation that~$Q_B$ is weakly connected and~$n\ge 2$.
Relabeling vertices $1,\dots,n$ if necessary, we may assume that
\begin{gather*}
\text{if} \ i\to j \ \text{is in} \ Q_B, \ \text{then} \ i<j.
\end{gather*}

Let $\{e_i\}$ be the standard basis of $\mathbb{Z}^n$, and $(r_1,\dots,r_n)$ is a permutation of $\{1,\dots,n\}$. We say that a map $f\colon \mathbb{Z}^n\to \mathbb{Z}^n$ is {\it coordinate-wise $\mathbb{Z}_{\ge0}$-linear} if
\begin{gather*}
f\left(\sum_i m_i\epsilon_ie_i\right)=\sum_i m_i f(\epsilon_ie_i),\quad \text{for all} \ \ m_i\in\mathbb{Z}_{\ge0},\ \  \epsilon_i\in\{1,-1\}.
\end{gather*}
We call such a map {\it triangularizable} if there is an order on the set $\{1,\dots,n\}$ as $\{r_1,\dots,r_n\}$ such that
\begin{gather*}
f(\epsilon e_{r_i})+\epsilon e_{r_i}\in \mathbb{Z}e_{r_{i+1}}\oplus\cdots \oplus\mathbb{Z}e_{r_n},\qquad \forall\, 1\le i\le n, \quad \epsilon\in\{1,-1\}.
\end{gather*}

\begin{Lemma}\label{lem:triangularizable bijective}
Let $f$ be a triangularizable coordinate-wise $\mathbb{Z}_{\ge0}$-linear map. Then $f$ is bijective.
\end{Lemma}

\begin{proof}
Without loss of generality we assume $r_i=i$ for $i=1,\dots,n$. We show that $f$ is invertible, i.e., for any $\sum c_ie_i\in\mathbb{Z}^n$,  there exists a unique $\sum b_ie_i\in\mathbb{Z}^n$ such that
\begin{equation}\label{fb=c}
f\left(\sum b_ie_i\right)=\sum c_ie_i.
\end{equation}

We f\/irst show that $f$ is injective. Assume $\sum b_ie_i$ satisf\/ies \eqref{fb=c}. Def\/ine a function
\begin{gather*}
s\colon \ \mathbb{Z}\to\{1,-1\}, \qquad  s(b)=
\begin{cases} \hphantom{-}1,&  \text{if} \ \ b\ge0,\\
  -1, & \text{otherwise},
\end{cases}
\end{gather*}
Denote the usual inner product on $\mathbb{Z}^n$ where $\langle e_i,e_j\rangle=1$ if $i=j$, otherwise $0$. Then for each $1\le k\le n$, we must have
\begin{gather*}
c_k =\left\langle\sum c_ie_i,e_k\right\rangle=\sum_{i=1}^{n}|b_i| \langle f(s(b_i) e_i),e_k\rangle
\\
\hphantom{c_k}{}
=\left(\sum_{i=1}^{k-1}|b_i| \langle f(s(b_i) e_i),e_k\rangle\right)+|b_k|\langle f(s(b_k) e_k),e_k\rangle=\left(\sum_{i=1}^{k-1}|b_i| \langle f(s(b_i) e_i),e_k\rangle\right)-b_k,
\end{gather*}
where the third equality is because of the following: for $i>k$, the triangularizable property assures that $f(\pm e_i)$ would be written only using a linear combination of $e_j$'s with $j>k$; since $\langle e_j,e_k\rangle=0$ for these $j$'s, we have $\langle f(\pm e_i),e_k\rangle=0$.
Therefore $b_1,\dots,b_n$ are uniquely determined (in that order) by the equation
\begin{gather}\label{bk}
b_k=\left(\sum_{i=1}^{k-1}|b_i| \langle f(s(b_i) e_i),e_k\rangle\right)-c_k,
\end{gather}
 implying that $f$ is injective.

To show that $f$ is surjective, we construct $b_1,\dots,b_n$, in that order, by~\eqref{bk}. Then we claim that~\eqref{fb=c} holds. Indeed,
\begin{gather*}
\left\langle f\left(\sum b_ie_i\right),e_k\right\rangle=\sum_{i=1}^n |b_i|\langle f(s(b_i) e_i),e_k\rangle=\left(\sum_{i=1}^{k-1}|b_i| \langle f(s(b_i) e_i),e_k\rangle\right)-b_k=c_k,
\end{gather*}
for all $1\le k\le n$, which implies~\eqref{fb=c}.
\end{proof}

Take the lexicographic order
\begin{gather*}
x_1>x_2>\cdots >x_n>1,
\end{gather*}
 that is, $\prod x_i^{a_i}< \prod x_i^{b_i}$ if the f\/irst pair of unequal exponents $a_i$ and $b_i$ satisfy $a_i<b_i$.

\begin{Lemma}\label{lem:l injective}
The map $f\colon \mathbb{Z}^n\to\mathbb{Z}^n$ sending ${\bf a}$ to the exponent vector of the lowest Laurent monomial of $z[{\bf a}]$ is triangularizable coordinate-wise $\mathbb{Z}_{\ge0}$-linear, so is injective.

Moreover, the coefficient of the lowest Laurent monomial of $z[{\bf a}]$ $($which a priori is a Laurent polynomial in $\mathbb{ZP})$ is a Laurent monomial in $\mathbb{ZP}$, i.e., invertible in $\mathbb{ZP}$.
\end{Lemma}

\begin{proof}
It follows from the def\/inition of $z[{\bf a}]$ that $f$ is coordinate-wise $\mathbb{Z}_{\ge0}$-linear.
Next we show that $f$ is triangularizable. Take $1\le k\le n$. By the def\/inition of mutation,
\begin{gather*}
x_k' = x_k^{-1} \Bigg( \underbrace{\prod_{1\leq i<k} x_i^{b_{ik}}}_{P_1} \underbrace{\prod_{i>n} x_i^{[b_{ik}]_+}}_{Q_1} + \underbrace{\prod_{k < i \leq n} x_i^{-b_{ik}}}_{P_2} \underbrace{\prod_{i>n} x_i^{[-b_{ik}]_+}}_{Q_2} \Bigg),
\end{gather*}
so its lowest monomial is
\begin{gather*}
 \begin{cases} x_k^{-1} P_2 Q_2, & \text{if $k$ is not a source in $Q_B$,}  \\ x_k^{-1} P_1Q_1=
x_k^{-1} Q_1, & \text{if $k$ is a source in $Q_B$.} \end{cases}
\end{gather*}
So
\begin{gather*}
 f(e_k)=\begin{cases} \displaystyle -e_k+\sum_{k<i\le n}(-b_{ik})e_i, & \text{if $k$ is not a source in $Q_B$,}  \\ -e_k, & \text{if $k$ is a source in $Q_B$.} \end{cases}
\end{gather*}
We also have the obvious equality $f(-e_k)=e_k$. Thus $f$ is triangularizable by taking the order $(1,2,\dots,n)$.
Therefore, $f$ is injective by Lemma~\ref{lem:triangularizable bijective}.

For the ``moreover'' statement, it suf\/f\/ices to note that $P_1$ and $P_2$ in the above def\/inition of~$x'_k$ have dif\/ferent  exponent vectors. This is because of the assumption that $n\ge 2$ and~$Q_B$ is weakly connected, so it is impossible for $P_1=P_2=1$ to hold.
\end{proof}

\begin{Lemma}\label{lem:lowest equal}
 For any ${\bf a} \in \mathbb{Z}^n$, the lowest Laurent monomials in $z[{\bf a}]$ and $\tilde{x}[{\bf a}]$ are equal. As a~consequence, $\{\tilde{x}[{\bf a}]\}_{{\bf a}\in\mathbb{Z}^n}$ is $\mathbb{ZP}$-linearly independent.
\end{Lemma}

\begin{proof}
First note that the consequence follows from Lemma~\ref{lem:l injective}, so for the rest we focus on the f\/irst statement.
By the multiplicative property of $z[{\bf a}]$ and $\tilde{x}[{\bf a}]$ (Lemma~\ref{multiplicative}), if suf\/f\/ices to consider the case when $\mathbf{a}=(a_i)\in \{0,1\}^n$. In the following, we use the notation from Lemma~\ref{lemma:01}. Def\/ine
\begin{gather*}
J:=\{1\le k\le n \,|\, a_k=1\},\qquad H:=\{1\le k\le n \,|\,  k\textrm{ is a source in } Q_B\}.
\end{gather*}
Since $z[{\bf a}]=\prod\limits_{k\in J} x'_k$, its lowest Laurent monomial is
\begin{gather}\label{lowest}
\left(\prod_{l=1}^n  x_l^{-a_l}\right)\left(\prod_{k\in J\setminus H} \left(\prod_{k < i \leq n} x_i^{-b_{ik}}\prod_{i>n} x_i^{[-b_{ik}]_+}\right) \right)\left( \prod_{k\in J\cap H} \prod_{i>n} x_i^{[b_{ik}]_+}\right).
\end{gather}
Def\/ine ${\bf s}=(s_1,\dots,s_n)$ by $s_k=1$ if $k\in J\setminus H$, otherwise take $s_k=0$.
 For $i\le n$, the degree of~$x_i$ in~\eqref{lowest} is
\begin{gather*}
-a_i+\sum_{k\in J\setminus H}[-b_{ik}]_+
=-a_i+\sum_{k=1}^ns_k[-b_{ik}]_+
=-a_i+
\sum_{k=1}^n \big(s_k[-b_{ik}]_++(a_k-s_k)[b_{ik}]_+\big),
\end{gather*}
where the second equality is because  $(a_k-s_k)[b_{ik}]_+=0$. Indeed, if $a_k-s_k\neq0$, then $a_k=1$ and $s_k=0$, thus $k\in J\cap H$. Then $k\in H$ implies $b_{ik}\le 0$,  therefore $[-b_{ik}]_+=0$.
For $i> n$, the degree of $x_i$ in \eqref{lowest} is
\begin{gather*}
\sum_{k\in J\setminus H}[-b_{ik}]_++\sum_{k\in J\cap H}[b_{ik}]_+
=\sum_{k=1}^ns_k[-b_{ik}]_++(a_k-s_k) [b_{ik}]_+.
\end{gather*}
So we can rewrite \eqref{lowest} as
\begin{gather*}
\left(\prod_{i=1}^n  x_i^{-a_i}\right)\prod_{i=1}^m
x_i^{\sum\limits_{j=1}^n (a_j-s_j)[b_{ij}]_+ + s_j[-b_{ij}]_+}.
\end{gather*}
It is easy to check that ${\bf s}$ is in ${\bf S}$. Comparing with~\eqref{df:x[a]} in Lemma~\ref{lemma:01} we draw the expected conclusion.
\end{proof}

\begin{Lemma}\label{supportofS}
Let $t\in\{1,\dots,n\}$ and  $M\in\mathbb{Z}$.
Assume that all the Laurent monomials $cx_1^{r_1}\cdots x_n^{r_n}$ $(c\neq0\in\mathbb{ZP})$ appearing in the sum $S:=\sum\limits_{\bf b} u({\bf b}) z[{\bf b}] $ $(u({\bf b}) \in \mathbb{ZP})$ satisfy $-r_t\le M$. Then $\langle{\bf b},e_t\rangle\le M$ for all ${\bf b}$ with $u({\bf b})\neq 0$ in~$S$.
\end{Lemma}

\begin{proof}
We f\/irst assume that $r_t\le 0$.
Take the lexicographic order
\begin{gather*}
x_t>x_{t-1}>\cdots >x_1>x_{t+1}>x_{t+2}>\cdots>x_n>1.
\end{gather*}
Consider the map  $f\colon \mathbb{Z}^n\to\mathbb{Z}^n$ that sends ${\bf b}$ to the  exponent vector of the lowest term of $z[{\bf b}]$.

We claim that $f$ is a triangularizable coordinate-wise $\mathbb{Z}_{\ge0}$-linear map by taking the order $(t,t-1,\dots,1,t+1,t+2,\dots,n)$. Indeed,  the proof is similar to the one of Lemma \ref{lem:l injective}, the only dif\/ference is that now we have{\samepage
\begin{gather*} f(e_k)=\begin{cases}
\displaystyle -e_k \ \ \text{or} \ \ -e_k+\sum_{i<k}b_{ik}e_i, & \text{if $k<t$,} \\
\displaystyle -e_k \ \ \text{or} \ \  -e_k+\sum_{k<i\le n}(-b_{ik})e_i, & \text{if $k>t$,}  \\
\displaystyle -e_k+\sum_{i<k}b_{ik}e_i \ \ \text{or} \ \   -e_k+\sum_{k<i\le n}(-b_{ik})e_i , & \text{if $k=t$.} \\
\end{cases}
\end{gather*}
(It is unnecessary for our purpose to f\/igure out the exact value of $f(e_k)$.)}

\looseness=-1
Thus we conclude that $f$ is bijective by Lemma~\ref{lem:triangularizable bijective}.
Let~$z[{\bf b}']$ be the one that appears in~$S$ whose lowest Laurent monomial $cx_t^{-\langle{\bf b}',e_t\rangle}\prod\limits_{i\neq t} x_i^{r_i}$ is smallest. Note that $\langle{\bf b},e_t\rangle\le\langle{\bf b}',e_t\rangle$ for all~${\bf b}$ appearing in~$S$. By the injectivity of $f$, the lowest Laurent monomial $cx_t^{-\langle{\bf b}',e_t\rangle}\prod\limits_{i\neq t} x_i^{r_i}$ will not  cancel out with other terms in~$S$. So $\langle{\bf b}',e_t\rangle\le M$ by the assumption that all the Laurent monomials $cx_1^{r_1}\cdots x_n^{r_n}$ that appear in~$S$ satisfy $-r_t\le M$. Thus $\langle{\bf b},e_t\rangle\le M$ for all~${\bf b}$ appearing in~$S$.
\end{proof}

Def\/ine a partial order $\prec$ on $\mathbb{Z}^n$:
\begin{gather*}
{\bf b} \prec {\bf a} \  \Leftrightarrow\  b_i\le a_i  \ \ \textrm{for all}  \ 1\le i\le n, \qquad \text{and} \qquad \sum_i[b_i]_+<\sum_i [a_i]_+.
\end{gather*}

\begin{Lemma}\label{lem:xinz}
For each ${\bf a} \in \mathbb{Z}^n$, we have
\begin{gather*}
\tilde{x}[{\bf a}] - z[{\bf a}] = \sum_{{\bf b} \prec {\bf a}} u({\bf a},{\bf b}) z[{\bf b}], \qquad
\text{where} \quad u({\bf a},{\bf b}) \in \mathbb{ZP}.
\end{gather*}
\end{Lemma}

\begin{proof}
Since $z[{\bf a}]$ form a $\mathbb{ZP}$-basis of the cluster algebra  $\mathcal{A}$, we can write
\begin{gather*}
\tilde{x}[{\bf a}] = \sum_{{\bf b}} u({\bf a},{\bf b}) z[{\bf b}],
\end{gather*}
where $u({\bf a},{\bf b})\neq 0\in\mathbb{ZP}$.
For $1\le t\le n$, it is easy to see that all the Laurent monomials $cx_1^{r_1}\cdots x_n^{r_n}$ ($c\neq0\in\mathbb{ZP}$) appeared $\tilde{x}[{\bf a}]$ satisfy $-r_t\le a_t$, thus  $b_t\le a_t$ by Lemma \ref{supportofS}.
Together with Lemma~\ref{lem:lowest equal} we can conclude that
\begin{gather}\label{weaker expansion x-z}
\tilde{x}[{\bf a}] - z[{\bf a}] = \sum_{\bf b} u({\bf a},{\bf b}) z[{\bf b}], \qquad u({\bf a},{\bf b}) \in \mathbb{ZP},
\end{gather}
where ${\bf b}\neq{\bf a}$ satisf\/ies $b_i\le a_i$ for all $1\le i\le n$.

In the rest we show that $\sum\limits_i[b_i]_+<\sum\limits_i [a_i]_+$. For simplicity, we assume that the f\/irst $n'$ coordinates of ${\bf a}$ are positive and the rest are non-positive. If $n'=n$, then $\sum\limits_i[b_i]_+<\sum\limits_i [a_i]_+$ since ${\bf b}\neq{\bf a}$. Assume $n'<n$. Regard $x_i$ for $n'<i\le n$ as frozen variable. Dividing~\eqref{weaker expansion x-z} by the monomial $x_{n'+1}^{-a_{n'+1}}\cdots x_{n}^{-a_{n}}$, we get the expansion of \begin{gather*}
\tilde{x}[(a_1,\dots,a_{n'})] - z[(a_1,\dots,a_{n'})]
 \end{gather*}
 in the basis $\{z[{\bf b}']\}_{{\bf b}\in\mathbb{Z}^{n'}}$ with coef\/f\/icients in $\mathbb{ZP}[x_{n'+1}^{\pm1},\dots, x_n^{\pm1}]=\mathbb{Z}[x_{n'+1}^{\pm1},\dots, x_m^{\pm1}]$. Using induction on $n$, we conclude that ${\bf b}'\neq(a_1,\dots,a_{n'})$, thus $\sum\limits_i[b_i]_+<\sum\limits_i [a_i]_+$.
\end{proof}

\begin{proof}[Proof of Theorem~\ref{thm:basis}]
The $\mathbb{ZP}$-independency of $\{\tilde{x}[{\bf a}]\}_{{\bf a}\in\mathbb{Z}^n}$ is asserted in Lemma~\ref{lem:lowest equal}.
Now we show that the $\mathbb{ZP}$-linear span of~$\{\tilde{x}[{\bf a}]\}_{{\bf a}\in\mathbb{Z}^n}$ equals~$\mathcal{A}$. Since $\{z[{\bf a}]\}_{{\bf a}\in\mathbb{Z}^n}$ is a~$\mathbb{ZP}$-basis of~$\mathcal{A}$, it suf\/f\/ices to show that, for any~${\bf a}\in\mathbb{Z}^n$, $z[{\bf a}]$~can be expressed as a $\mathbb{ZP}$-linear combination of~$\tilde{x}[{\bf a}']$ for ${\bf a'}\in\mathbb{Z}^n$.  Iteratively applying Lemma~\ref{lem:xinz}, and noting that~$z[{\bf b}]=\tilde{x}[{\bf b}]$ for~${\bf b}\in\mathbb{Z}^n_{\le0}$ (actually the equality holds if ${\bf b}$ has at most one positive coordinate),  we can express~$z[{\bf a}]$ as a~linear combination of~$\tilde{x}[{\bf b}]$'s.
\end{proof}

\subsection{The second proof}\label{section6.2}

Def\/ine a partial order on $\mathbb{Z}^n  $ by
\begin{gather}\label{eq:prec2}
{\bf b} \prec {\bf a} \ \ \text{if and only if}\ \ \sum_{i=1}^n [b_i]_+ < \sum_{i=1}^n [a_i]_+.
\end{gather}

\begin{Definition}
For ${\bf a}\in\mathbb{Z}^n$, let $T^{\bf a}$ be the set of all quasi-compatible collections and let $T^{\bf a}_{\rm gcc} \subset T^{\bf a}$ be the globally compatible collections. (See Def\/inition~\ref{definition:GCC} for the terminology.)
\end{Definition}

Let $ {\bf a} \in \{0,1\}^n  $. Let $C^{(i,j)} = \{ u_1^{(i,j)} , v_1^{(i,j)} \}$ be the corner associated to the edge $(i,j)$ (if it exists). Let $\overline{T}^{\bf a}=T^{\bf a}\setminus T^{\bf a}_{\rm gcc}$ and let $t$ be a quasi-compatible collection in $\overline{T}^{\bf a}$. Let $c(t)$ be the set of all corners of $t$.
 We can def\/ine an equivalence relation on quasi-compatible collections $t,s \in \overline{T}^{\bf a}$ by $t \sim s$ if and only if $c(t)=c(s)$. Let $M_t \subset \overline{T}^{\bf a}$ denote the equivalence class of $t$ with respect to $\sim.$
 To each equivalence class we can associate a vector
 \begin{gather*}
 {\bf c}_t = \sum_{C^{(i,j)} \in c(t)}\left(\sum_{(i,h)\in Q_B} b_{ih} e_h + \sum_{(k,j)\in Q_B} (-b_{jk}) e_k\right) \in \mathbb{Z}_{\geq 0}^n.
 \end{gather*}

For each $M_t \subset \overline{T}^{\bf a}$ we also associate a vector ${\bf b}_t={\bf a}-{\bf c}_t$. Note that ${\bf b}_t \prec {\bf a}$. Indeed, since $t\notin T^{\bf a}_{\rm gcc}$, there is at least one corner $C^{(i,j)}$ in $c(t)$,
then $({\bf b}_t)_i\le a_i-(-b_{ji})<a_i=1$, $({\bf b}_t)_j\le a_j-b_{ij}<a_j=1$, and
 $({\bf b}_t)_k\le a_k$ for all $k$, thus $\sum\limits_{r=1}^n [({\bf b}_t)_r]_+ < \sum\limits_{r=1}^n [a_r]_+$.

Let $\mathcal{D}_{\bf a}$ be the set of edges of Dyck paths for the vector ${\bf a}$. Let $\mathcal{D}_{{\bf b}_t}$ be def\/ined analogously. Since ${\bf b}_t \prec {\bf a}$, $\mathcal{D}_{{\bf b}_t} \subset \mathcal{D}_{\bf a}$ (since each horizontal (resp.~vertical) edge in  $\mathcal{D}_{{\bf b}_t}$ naturally corresponds to  a horizontal (resp. vertical) edge in $\mathcal{D}_{\bf a}$), and by the def\/inition of ${\bf b}_t$ we have
\begin{gather*}
\mathcal{D}_{\bf a} \setminus \mathcal{D}_{{\bf b}_t} =  \bigcup_{C^{(i,j)}\subset s \in M_t}\big\{ v_1^{(k,j)},u_1^{(i,h)},u_1^{(j,f)},v_1^{(g,i)} \,|\, (j,f),(g,i),(i,h),(k,j) \in Q_{\tilde B} \big\}.
\end{gather*}

\begin{Definition} Def\/ine the map $\phi_{{\bf b}_t}\colon T^{{\bf b}_t} \rightarrow \overline{T}^{\bf a}$ as follows
\begin{gather*}
\phi_{{\bf b}_t}(d)= \left( \bigcup_{C^{(i,j)}\subset s \in M_t}\big\{u_1^{(i,h)}, v_1^{(k,j)} \,|\, (i,h),(k,j) \in Q_{\tilde B} \big\}\right) \sqcup d.
\end{gather*}
\end{Definition}

It is easy to see that if $d \in T^{{\bf b}_t}$ is quasi-compatible then so is $\phi_{{\bf b}_t}(d) \in \overline{T}^{\bf a}$ and that $\phi_{{\bf b}_t}$ is injective.

\begin{Lemma}\label{lem:phi}
Let ${\bf a} \in \{0,1\}^n$, and $t$, $M_t$, ${\bf b}_t$, $\phi_{{\bf b}_t}$ be defined as above.
\begin{enumerate}\itemsep=0pt

\item[$1.$]  If there exists a quasi-compatible collection $d\in T^{{\bf b}_t}$ such that $\phi_{{\bf b}_t}(d) \in M_t$, then $M_t$ is contained in the image of $\phi_{{\bf b}_t}$. Furthermore, $ \{\Ima \phi_{{\bf b}_t} |M_t\subset \bar{T}^{\bf a}\}$ cover $\overline{T}^{\bf a}$.

\item[$2.$] For $A \subset T^{\bf a}$ define
\begin{gather*}
z[A] := \left(\prod_{i=1}^n  x_i^{-a_i}\right)
\sum _{s \in A}\left(\prod_{(i,j)\in Q_{\tilde B}}  x_{i}^{b_{ij} \big|S^{(i,j)}_2 \big|} x_{j}^{-b_{ji} \big|S^{(i,j)}_1 \big|}\right).
\end{gather*}
Then $z[\Ima \phi_{{\bf b}_t}]=z[{\bf b}_t]$. It follows immediately that if $\Ima \phi_{{\bf b}_t}=M_t$ then~$z[M_t]$ is a standard monomial.
\item[$3.$] We can define a partial order on~$\{ M_t\}$ by inclusion of sets of corners. If~$M_t$ is maximal with respect to this order then $\Ima \phi_{{\bf b}_t}=M_t$.
\end{enumerate}
\end{Lemma}

\begin{proof}
(1)~To see that the $\{\Ima \phi_{{\bf b}_t}\}$ cover we will show that $M_t \subset \Ima\phi_{{\bf b}_t}$.  It suf\/f\/ices to f\/ind a~single quasi-compatible collection that maps into~$M_t$, because
If $\phi_{{\bf b}_t}(d) \in M_s$ then any~$w$ in~$M_s$ is a disjoint union of
\begin{gather*}
\bigcup_{C^{(i,j)}\subset s \in M_t}\big\{ v_1^{(k,j)},u_1^{(i,h)} \, |\, (i,h),(k,j) \in Q_{\tilde B} \big\}
\end{gather*}
 with an element in $T^{{\bf b}_t}$.
To produce this collection take a globally compatible collection in~$T^{{\bf b}_t}$ that contains no~$u_1^{(k,j)}$ or~$v_1^{(i,h)}$ for all~$i$,~$j$ such that $C^{(i,j)} \subset t$ and for all $(k,j),(i,h) \in Q_{\tilde B}$.  This property guarantees that no new corners  besides those in~${\bf c}(t)$ will be created when we map it into~$\overline{T}^{\bf a}$. Therefore its image under $\phi_{{\bf b}_t}$ is in~$M_t$.

(2)~Using the fact that every quasi-compatible collection in the image of $\phi_{{\bf b}_t}$ shares the edges
\begin{gather*}
\bigcup_{C^{(i,j)}\subset s \in M_t}\big\{ v_1^{(k,j)},u_1^{(i,h)} \,|\,(i,h),(k,j) \in Q_{\tilde B} \big\},
\end{gather*}
 and the def\/inition of~${\bf b}_t$ it follows that
\begin{gather*}
z[{\bf b}_t]  = \left(\prod_{i=1}^n  x_i^{-({\bf b}_t)_i}\right)
\sum_{s \in T^{{\bf b}_t}}\left(\prod_{(i,j)\in Q_{\tilde B} }  x_{i}^{b_{ij}\big|S^{(i,j)}_2\big|} x_{j}^{-b_{ji}\big|S^{(i,j)}_1\big|}\right) \\
\hphantom{z[{\bf b}_t]}{}
 =  \left(\prod_{i=1}^n  x_i^{-a_i}\right)\left(\prod_{i=1}^n  x_i^{({\bf c}_t)_i}\right)
\sum_{s \in T^{{\bf b}_t}}\left(\prod_{(i,j)\in Q_{\tilde B} }  x_{i}^{b_{ij}\big|S^{(i,j)}_2\big|} x_{j}^{-b_{ji}\big|S^{(i,j)}_1\big|}\right) \\
\hphantom{z[{\bf b}_t]}{}
= \left(\prod_{i=1}^n  x_i^{-a_i}\right)
\sum_{s \in \Ima \phi_{{\bf b}_t}}\left(\prod_{(i,j)\in Q_{\tilde B} }  x_{i}^{b_{ij}\big|S^{(i,j)}_2\big|} x_{j}^{-b_{ji}\big|S^{(i,j)}_1\big|}\right)
 = z[\Ima \phi_{{\bf b}_t}].
\end{gather*}

(3) We have already shown that $M_t \subset \Ima \phi_{{\bf b}_t}$. Since $M_t$ is maximal, ${\bf c}(t)$~is a maximal set of corners and every element in the image of $\phi_{{\bf b}_t}$ contains ${\bf c}(t)$. So $ \Ima \phi_{{\bf b}_t}\subset M_t.$
\end{proof}

We need the following facts about $z[{\bf a}]$.
\begin{Lemma}\label{zproduct}
For any ${\bf a},{\bf b}\in\mathbb{Z}^n$,
\begin{gather*}
z[{\bf a}]z[{\bf b}]=\sum_{\bf c} u({\bf c})z[{\bf c}],\quad u({\bf c})\in\mathbb{ZP},
\end{gather*}
where ${\bf c}$ satisfies $[c_t]_+\le \sum[a_t+b_t]_+$ for $1\le t\le n$.
\end{Lemma}

\begin{proof}
Regard $x_i$ ($i\neq t$) as frozen variables, i.e., $x_t$ as the only non-frozen variable. Then it is reduced to proving the lemma for $n=1$, which is easy.
\end{proof}

Now we can prove the following lemma.
\begin{Lemma}\label{lem:xinz2}
For ${\bf a} \in \mathbb{Z}^n$  the expansion of $\tilde{x}[{\bf a}]$ in the basis of standard monomials is of the form
\begin{gather*}
 \tilde{x}[{\bf a}] - z[{\bf a}] = \sum_{{\bf b} \prec {\bf a}} u({\bf a},{\bf b}) z[{\bf b}],
 \end{gather*}
where $u({\bf a},{\bf b}) \in \mathbb{ZP}$ and all but finitely many are zero, and~$\prec$ as is defined in~\eqref{eq:prec2}.
\end{Lemma}

\begin{proof}
\emph{Step~I.} We prove the case when $m=n$, i.e., $\mathbb{ZP}=\mathbb{Z}$.

By the multiplicative property of $\{\tilde{x}[{\bf a}]\}$ and $\{z[{\bf a}]\}$ (Lemma~\ref{multiplicative}) it is suf\/f\/icient to show the result for ${\bf a}\in \{0,1\}^n$.
Indeed, if ${\bf a}={\bf a}'+{\bf a}''$ and we have the result for ${\bf a}'$, and ${\bf a}''$ then{\samepage
\begin{gather*}
\tilde{x}[{\bf a}] =\tilde{x}[{\bf a}']\tilde{x}[{\bf a}'']
 =\left(z[{\bf a}']+\sum_{{\bf b}'\prec {\bf a}'} u({\bf a}',{\bf b}') z[{\bf b}']\right)\left(z[{\bf a}'']+\sum_{{\bf b}''\prec {\bf a}''} u({\bf a}'',{\bf b}'') z[{\bf b}'']\right)\\
\hphantom{\tilde{x}[{\bf a}]}{}
\stackrel{(*)}{=}z[{\bf a}']z[{\bf a}'']  +\sum_{{\bf b}\prec {\bf a}'+{\bf a}''} u({\bf b}) z[{\bf b}]
 = z[{\bf a}] +\sum_{{\bf b}\prec {\bf a}} u({\bf b}) z[{\bf b}],
\end{gather*}
where $(*)$ uses Lemma \ref{zproduct}.}

We have that
\begin{gather*}
z[{\bf a}]=z[T^{\bf a}]=z[T^{\bf a}_{\rm gcc}]+z[\overline{T}^{\bf a}]=\tilde{x}[{\bf a}]+z[\overline{T}^{\bf a}].
\end{gather*}
So we only need to write $z[\overline{T}^{\bf a}]$ in terms of standard monomials.
Now $\{M_t\}$ forms a partition of~$\overline{T}^{\bf a}$. So,
\begin{gather*}
z[\overline{T}^{\bf a}] = \left(\prod_{i=1}^n  x_i^{-a_i}\right)
\sum_{s \in \overline{T}^{\bf a}} \left(\prod_{(i,j)\in Q_{\tilde B}}  x_{i}^{b_{ij}\big|S^{(i,j)}_2\big|} x_{j}^{-b_{ji}\big|S^{(i,j)}_1\big|}\right)   \\
\hphantom{z[\overline{T}^{\bf a}]}{}
= \left(\prod_{i=1}^n  x_i^{-a_i}\right) \sum_{M_t \subset \overline{T}^{\bf a}} \sum_{s \in M_t} \left(\prod_{(i,j)\in Q_{\tilde B}}  x_{i}^{b_{ij}\big|S^{(i,j)}_2\big|} x_{j}^{-b_{ji}\big|S^{(i,j)}_1\big|}\right)
 = \sum_{ M_t \subset \overline{T}^{\bf a}} z[M_t].
\end{gather*}
Note that the $z[M_t]$ are not necessarily standard monomials. However by Lemma~\ref{lem:phi}, $z[\Ima \phi_{{\bf b}_t}]$ are standard monomials, and $\{M_t\}$ partition the image of $\phi_{{\bf b}_t}$. Therefore we can do the same manipulation as above to get the equation
\begin{gather}\label{zbt}
z[{\bf b}_t] = z[\Ima \phi_{{\bf b}_t}]=\sum_{ M_s \subset \Ima \phi_{{\bf b}_t}} z[M_s] .
\end{gather}
Now there are f\/initely many ${\bf b}_t$ and ${\bf b}_t \prec a$ (since~${\bf c}_t \in \mathbb{Z}^n_{\geq 0}$), so to f\/inish the proof of our claim we only need to show that~$z[\overline{T}^a]$ is a~$\mathbb{Z}$-linear combination of  these~$z[{\bf b}_t]$'s.

Note that $M_s \subset \Ima \phi_{{\bf b}_t}$ has at least the corners of $c(t)$ but could have more. So $M_s \not\subset \Ima\phi_{b_t}$ if $M_s < M_t$ with respect to our order in Lemma~\ref{lem:phi}(3), and we can rewrite~\eqref{zbt} as
\begin{gather*}
z[b_t]=z[M_t]+\sum_{M_t<M_s} d_sz[M_s], \qquad d_s\in\{0,1\}.
 \end{gather*}
 By extending the partial order on~$\{M_t\}$ to a total order, we see that the transition matrix between~$\{z[{\bf b}_t]\}$ and~$\{z[M_t]\}$ is triangular with all diagonal entries~1, so it is invertible over~$\mathbb{Z}$.
Therefore every~$z[M_t]$ is a linear combination of~$z[{\bf b}_t]$.

Therefore, $z[\overline{T}^{\bf a}]$ is a $\mathbb{Z}$-linear combination of~$z[{\bf b}_t]$, where~$t$ satisf\/ies $M_t \in \overline{T}^{\bf a}$.

\emph{Step II.} We prove the principal coefficient case, i.e., $m=2n$ and $\tilde{B}=\begin{bmatrix}
B\\ I_n
\end{bmatrix}$ where $I_n$ is the $n\times n$ identity matrix.

Let   $\mathbb{ZP}=\mathbb{Z}[y_1,\dots,y_n]$ where $y_j=x_{n+j}$ for $1\le j\le n$.
Note that this will imply the general coef\/f\/icient case by replacing the principal coef\/f\/icients $y_j$ by $\prod\limits_{i=n+1}^mx_i^{b_{ij}}$ for $1\le j\le n$.

Apply Step I to the coef\/f\/icient-free cluster algebra $\mathcal{A}'$ with $B$-matrix
$\begin{bmatrix}
B&-I_n\\
I_n&0
\end{bmatrix}$. For any ${\bf a}\in\mathbb{Z}^n$, denote $\tilde{\bf a}=(a_1,\dots,a_n,0,\dots,0)\in\mathbb{Z}^{2n}$. Then $\tilde{x}[\tilde{\bf a}]=\tilde{x}[{\bf a}]$ and $z[\tilde{\bf a}]=z[{\bf a}]$. Thus
\begin{gather*}
 \tilde{x}[{\bf a}] - z[{\bf a}]=\tilde{x}[\tilde{\bf a}] - z[\tilde{\bf a}] = \sum u({\bf b}') z[{\bf b}'],
 \end{gather*}
where $u({\bf b}') \in \mathbb{Z}$ and ${\bf b}'=(b_1,\dots,b_{2n})$ satisf\/ies $\sum\limits_{i=1}^m [b_i]_+<\sum\limits_{i=1}^n[a_i]_+$. So it suf\/f\/ices to show that $z[{\bf b}']$ is a $\mathbb{ZP}$-linear combination of $z[{\bf c}]$ with ${\bf c}=(c_1,\dots,c_n)\in \mathbb{Z}^n$ satisfying $\sum\limits_{i=1}^n[c_i]_+\le \sum\limits_{i=1}^n[b_i]_+$ (so $\sum\limits_{i=1}^n[c_i]_+\le \sum\limits_{i=1}^m[b_i]_+<\sum\limits_{i=1}^n[a_i]_+$).

Denote ${\bf b}=(b_1,\dots,b_n)$. Since
$z[{\bf b}']=z[{\bf b}] \, z[(0,\dots,0,b_{n+1},\dots,b_{2n})]$ where the second factor is a $\mathbb{ZP}$-linear combination of $z[{\bf d}]$ with ${\bf d}=(d_1,\dots,d_n)\in \mathbb{Z}_{\le 0}^n$, it suf\/f\/ices to show that $z[{\bf b}]z[{\bf d}]$ is a $\mathbb{ZP}$-linear combination of  $z[{\bf c}]$ with ${\bf c}=(c_1,\dots,c_n)\in \mathbb{Z}^n$ satisfying $\sum[c_i]_+\le \sum[b_i]_+$. This follows from Lemma~\ref{zproduct} since $\sum\limits_{i=1}^n[c_i]_+ \le \sum\limits_{i=1}^n [b_i+d_i]_+\le \sum\limits_{i=1}^n [b_i]_+$.
\end{proof}

\begin{proof}[Proof of Theorem~\ref{thm:basis}]
Similar to the proof in Section~\ref{section6.1}, except that we use  Lemma~\ref{lem:xinz2} in place of Lemma~\ref{lem:xinz}.
\end{proof}

\subsection*{Acknowledgements}

The authors are grateful to the anonymous referees for carefully reading through the manuscript and giving us many constructive suggestions to improve the presentation. KL is supported by  Wayne State University, Korea Institute for Advanced Study, AMS Centennial Fellowship and NSA grant H98230-14-1-0323.
MM is supported by GAANN Fellowship.


\pdfbookmark[1]{References}{ref}
\LastPageEnding

\end{document}